\documentclass{jgcc} %%% last changed 2020-01-17

\pdfoutput=1
\usepackage{lastpage}
\jgccdoi{12}{1}{3} 
\jgccheading{}{\pageref{LastPage}}{}{}%
{Oct.~25,~2018}{Apr.~15,~2020}{}


\usepackage{hyperref}
\usepackage{amssymb}
\usepackage{amsthm}
\usepackage{calc}
\usepackage[centertags]{amsmath}
\usepackage{latexsym}
\usepackage{amsfonts}
\usepackage{multirow}
\usepackage{color}

\newcommand{\Z}{\mathbb{Z}}
\newcommand{\Q}{\mathbb{Q}}

\newcommand{\F}{\mathbb{F}}
\newcommand{\G}{\mathbb{G}}
\let\phi\varphi

\begin{document}

\title[Computing Optimal Ate Pairings]{Computing Optimal Ate Pairings on Elliptic Curves with Embedding Degree 9, 15 and 27}
%\titlecomment{{\lsuper*}}

\author[E.~Fouotsa]{Emmanuel Fouotsa}	%required
\address{Department of Mathematics, Higher Teacher Training College, The University of Bamenda P.O Box 39 Bambili, Cameroon}	%required
\email{emmanuelfouotsa@yahoo.fr}  %optional
\thanks{The first and third authors acknowledge support from The Simons Foundations through Pole of Research in Mathematics with applications to Information Security, Sub-Saharan Africa.}	%optional

\author[N.~El Mrabet]{Nadia El Mrabet}	%optional
\address{Mines Saint-Etienne, CEA-Tech, Centre CMP, Departement SAS, France}	%optional
\email{nadia.el-mrabet@emse.fr}  %optional
%\thanks{thanks 2, optional.}	%optional

\author[A.~Pecha]{Aminatou Pecha}	%optional
\address{Ecole Nationale Sup\'{e}rieur Polytechnique, Universit\'{e} de Maroua P.O.Box 46 Maroua, Cameroon}	%optional
%\urladdr{name3@url3\quad\rm{(optionally, a web-page can be specified)}}  %optional
\email{aminap2001@yahoo.fr}	%optional

%% etc.

%% required for running head on odd and even pages, use suitable
%% abbreviations in case of long titles and many authors:

%%%%%%%%%%%%%%%%%%%%%%%%%%%%%%%%%%%%%%%%%%%%%%%%%%%%%%%%%%%%%%%%%%%%%%%%%%%

%% the abstract has to PRECEDE the command \maketitle:
%% be sure not to issue the \maketitle command twice!

\begin{abstract}
%  \noindent
Much attention has been given to the efficient computation of pairings on elliptic curves with even embedding degree since the advent of pairing-based cryptography. The few existing works in the case of odd embedding degrees require some improvements.
  This paper considers the computation of optimal ate pairings on elliptic curves of embedding degrees $k=9$, $15$, $27$ which have twists of order three. Our main goal is to provide a detailed arithmetic and cost estimation of operations in the tower extensions field of the corresponding extension fields. A good selection of parameters 
  enables us to improve the theoretical cost for the Miller step and the final exponentiation using the lattice-based method as compared to the previous few works that exist in these cases. In particular, for $k=15$, $k=27$, we obtain an improvement, in terms of operations in the base field, of up to 25\% and 29\% respectively in the computation of the final exponentiation. 
  We also find that elliptic curves with embedding degree $k=15$ present faster results than BN12  curves at the 128-bit security level. 
  We provide a MAGMA implementation in each case to ensure the correctness of the formulas used in this work.
\end{abstract}

\maketitle

\section{Introduction}
Pairings are bilinear maps defined on the group of rational points of elliptic or hyper elliptic curves \cite{Was08}. They enable the realization of many cryptographic protocols such as the Identity-Based Cryptosystem \cite{BonFra01}, Identity-Based Encryption \cite{Coc01}, the Identity-Based Undeniable Signature \cite{LibQui04}, Short Signatures \cite{BonLynSha04} and Broadcast Encryption \cite{GovPanSahWat06}. Surveys of some applications of pairings can be found in \cite{DutBarSar04} and \cite[Chapter X]{BlaSerSma05}. These many applications justify the research on the efficient computation of pairings. Generally, let $E$ be an ordinary elliptic curve defined over a finite field $\F_p$ and $r$ be a large prime divisor of the order of the group $E(\F_p)$. The embedding degree of $E$ with respect to $r$ and the prime number $p$ is the smallest integer $k$ such that $r$ divides  $p^k-1$. The $\rho $-value of the elliptic curve $E$ is the value $\mbox{log}p/\mbox{log}r$ measuring the size of the base field relatively to the size $r$ of a subgroup of $E(\F_p)$. The Tate pairing and its variants are the most used pairings in cryptography. They map two linearly independent points of the subgroup of order $r$ of $E(\F_{p^k})$ to the group of $r$-th roots of unity in the finite field $\F_{p^k}$. The computation of the Tate pairing and its variants consists of an application of the Miller algorithm \cite{Mil04} and a final exponentiation. Efficient computation of pairings requires the construction of pairing-friendly elliptic curves over $\F_p$ with prescribed embedding degree $k$ (see for example \cite{BarNar05} or \cite{FreScoTes10}) and efficient arithmetic in the towering fields associated to $\F_{p^k}$ (see \cite{KobMen05}, \cite{GraSco10}, \cite{Kar13}, \cite{DevScoDah07}, \cite{FOUOBAN}). A lot of work has been done for shortening the Miller loop leading to the concept of pairing lattices \cite{Hes08},
%\cite{BarKimLynSco02}, and the latest is
or the optimal pairing described by Vercauteren which can be computed with the smallest number of iterations in the Miller algorithm \cite{Ver10}. Due to these advances, the final exponentiation step has became a serious task. In this work, we concentrate on elliptic curves $E$ over $\F_p$ with embedding degrees $9,15$ and $27$. 
\begin{table}[h]
	\centering
	\begin{tabular}{|c|c|c|c|c|}
		\hline
		% after \\: \hline or \cline{col1-col2} \cline{col3-col4} .
		Security & Bit length of  & Bit length of    &   $k$  &   $k$   \\
		level  &    $r$       &    $p^k$            &    $\rho\thickapprox 1$ &    $\rho\thickapprox 2$   \\
		%&$p^k$ & Bit length of\\$r$ & $k$, $\rho\thickapprox 1$& $k$, $\rho\thickapprox 2$\\
		\hline
		$80$ & $160$ & $960-1280$ & $6-8$ & $3-4$\\
		\hline
		$128$ & $256$ & $3000-5000$ & $12-20$& $6-10$ \\
		\hline
		$192$& $384$ & $8000-10000$ & $20-26$& $10-13$\\
		\hline
		$256$ & $512$& $14000-18000$ & $28-36$& $14-18$\\
		\hline
	\end{tabular}
	\caption{Bit sizes of curves parameters and corresponding embedding degrees to obtain commonly desired levels of security.}\label{recomendationParameters}
\end{table}
These curves admit twists of degree three which enable computations to be done in subfields and also lead to the denominator elimination technique. To our knowledge only few works (\cite{LeTan12}, \cite{MraGuiIon09} and \cite{ZhaLin12}) exist in these cases and much attention have only been given to elliptic curves with even embedding degree (see for example \cite{aranha12implementing},\cite{LauKnaRod11}, \cite{FOUOLOU}, \cite{FOUOBAN}). Another motivation for our work is the recent results on the resolution of the discrete logarithm problem \cite{KimBar16}. Indeed, according to the first analysis of this article, as for instance in \cite{Guille2016} and \cite{DuqBar18}, the security level for elliptic curves with friable embedding degree should be taken greater than those presented in Table \ref{recomendationParameters}. The main consequence is that elliptic curves with embedding degree $12$ or $18$ may not be the one ensuring a nice ratio between the security level and the arithmetic. Elliptic curves with odd embedding degree could become interesting and more efficient than elliptic curves with even embedding degree. We also  noticed that elliptic curves with odd embedding degree, especially with $k=27$, may be suitable for computing products of pairings \cite{ZhaLin12}.
In this work we consider  the following  parameter sizes  $(p\approx 2^{343},\,r\approx 2^{257},\, \rho=1.33)$,   $(p\approx 2^{575},\, r\approx 2^{385},\, \rho=1.5)$,  $(p\approx 2^{579},\, r\approx 2^{514},\, \rho=1.12)$ for curves with $k=9,\,15,\,27$ respectively. This corresponds to the 128, 192 and 256-bit security levels respectively according to the recommendations in Table \ref{recomendationParameters} \cite{FreScoTes10}. However, considering the recent recommendations based on the advances in Discrete Logarithm computation with the Number Field Sieve (NFS) algorithm and its variants (\cite{KimBar16}, \cite{RazMor15}, \cite{MenSarSha16}, \cite{DuqBar18}) the security level provided by the above parameters may reduce. Indeed, let $2^{-d}\mbox{exp}((c+o (1))(\mbox{log}Q)^{1/3}(\mbox{ log }\mbox{log}Q)^{2/3}$ where $d$ and $c$ are constants, be the running time of the NFS algorithm, with $Q=p^k$. The base-two logarithm of this runtime  with $(o(1)=0)$ gives  $S(Q,c,d)=c(\mbox{lg}e)(\mbox{log}Q)^{1/3}(\mbox{log}\mbox{log}Q)^{2/3}-d$. We then know that the constants $c$ and $d$, the embedding degree $k$ and the security level $l$ must satisfy the Pollard-Rho security and variants of NFS security constraints $\mbox{log}Q/k\ge 2\rho l$ and $S(Q,c,d)\ge l$.  Therefore the previous parameters provide a security level of 109, 168, 214-bits instead of 128, 192, 256-bit respectively for curves with $k=9$, $15$, $27$.
But for now we still consider in Sections \ref{k=9}, \ref{k=15} and  \ref{k=27} recommendation from Table \ref{recomendationParameters} in order to make a fair  comparison with previous work. Later, in Section \ref{update}, we consider advances in discrete logarithm computation and provide tentative updated parameters at the 128, 192, 256 bits level respectively for curves with $k=9$, $15$, $27$.
%and also in agreement with Menezes et \emph{al.} \cite{MenSarSha16} \emph{"...much work remains to be done before the impact of these new algorithms on concrete keylengths for pairing-based cryptography can be determined with full confidence...."}. 
So we propose a detailed arithmetic in the towering fields associated to the fields  $\F_{p^9},\F_{p^{15}}$ and $\F_{p^{27}}$. The lattice-based method explained by Fuentes \emph{et al.} \cite{LauKnaRod11} is applied to compute the final exponentiation in the cases $k=9$, 15. We also find a simple expression and explicit cost evaluation for the optimal pairing in the cases $k=9$ and $k=15$ as compared to the work in \cite{MraGuiIon09}. The results obtained are an improvement with respect to previous works \cite{LeTan12}, \cite{MraGuiIon09} and \cite{ZhaLin12} respectively for $k=9$, $15$ and $27$. Precisely, our contributions (see Table \ref{COSTCOMARISON} and Subsection \ref{compare} for comparison) in this work are:
\begin{enumerate}
	\item Determination of an explicit cost of the computation of the optimal pairing for the elliptic curves stated above. This includes a good selection of parameters for a shorter Miller loop and an efficient exponentiation. In particular, we saved one inversion in $\F_{p^{27}}$ for the computation of the Miller loop in the case $k=27$.
	\item Details of the arithmetic in the tower of subfields of $\F_{p^9},\F_{p^{15}}$ and $\F_{p^{27}}$.  In particular, we give the cost of the computation of Frobenius maps and Inversions in the cyclotomic subgroups of $\F_{p^9}^\ast,\F_{p^{15}}^\ast$ and $\F_{p^{27}}^\ast$, (see Appendices \ref{arithFP9}, \ref{arithFP27} and \ref{arithFP15}).
	\item Improvement in the costs of the final exponentiation by saving $828M_1+145S_1, 1170M_{1}+7767S_{1}$ and $8676M_{1}+32136S_1$ operations for elliptic curves of embedding degrees 9, 15 and 27 respectively, as compared to previous works in these cases; where $M_k$ and $S_k$ represent the costs of multiplication and squaring in the finite field $\F_{p^k}$.
	\item In Section \ref{update} we look for new  parameters considering the advances in Discrete Logarithm computation to update the cost of the optimal ate pairings on the studied curves at the 128, 192 and 256-bit security levels. We then compare our results with known curves such as BN and BLS curves. In particular we find that elliptic curves with embedding degree $k=15$ present faster results than BN12 curves at the 128-bit security levels. 
\end{enumerate}
We also provide a MAGMA implementation in each case to ensure the correctness of the formulas used in this work. The code is available in \cite{FouNadAmi18}.

The rest of this paper is organised as follows. In Section \ref{intro} we briefly present the Tate and ate pairings together with the Miller algorithm for their efficient computation; 
we also recall the concept of optimal ate pairing and the lattice-based method for computing the final exponentiation. Sections \ref{k=9}, \ref{k=15} and \ref{k=27} present arithmetic in subfields, and costs estimation of the Miller step and the final exponentiation when considering the embedding degrees $k=9$, $15$, $27$ respectively. Each of these sections includes a comparative analysis with previous work. 
Section \ref{fcomparaison} presents a general comparison of the results obtained in this work and the previous results in the literature. In Section \ref{update} we look for new  parameters considering the advances in Discrete Logarithm computation to update the cost of the optimal ate pairings on the studied curves at the 128, 192 and 256-bit security levels. We then compare our results with known curves such as BN and BLS curves.
We conclude the work in Section \ref{conclusion} in which we suggest as future work the search for parameters to have \textit{subgroup secure} ordinary curves \cite{BarCosMisNae15} and to ensure protection against \textit{small-subgroup attacks} \cite{LimLeeL97}.

Throughout the rest of this paper, we denote $M_k$, $S_k$ and $I_k$ respectively as the Cost of a multiplication, a squaring and an inversion in the field $\F_{p^k}$, for any integer $k$.

\section{Background and previous works}\label{intro}

%%%%%%%%%%%%%%%%%%%%%%
\subsection{Pairings and the Miller Algorithm}\label{reviewpairing}

Let $E$ be an elliptic curve defined over $\F_p$, a finite field of characteristic $p>3$. Let $r$ be a large prime factor of the group order of the elliptic curve. Let $m\in\Z$ and $P \in E(\F_{p})[r]$ a point of $E$ of order $r$ with coordinates in $\F_{p}$. Let $f_{m,P}$ be a function with divisor ${\rm Div}(f_{m,P})=m(P)-([m]P)-(m-1)(\mathcal{O})$ where $\mathcal{O}$ denotes the identity element of the group of points of the elliptic curve. Let $k$ denote the smallest integer such that $r$ divides $p^k-1$; this is called the embedding degree of $E$ with respect to $r$. We also consider the point  $Q\in E(\F_{p^{k}})[r]$ of $E$ of order $r$ with coordinates in $\F_{p^k}$ and let $\mu_{r}$  denote the group of $r$-th roots of unity in $\F^{\ast}_{p^{k}}$. The reduced Tate pairing $e_{r}$ is a bilinear and non-degenerate map defined as
\begin{align*}
e_r\colon  E(\F_p)[r] \times E(\F_{p^k})[r] &\longrightarrow  \mu_r, \\
(P,Q)  &\longmapsto  f_{r,P}(Q)^{\frac{p^k-1}{r}}.
\end{align*}
To define a variant of the Tate pairing called ate pairing \cite{HesSmaVer06}, we denote by $[i]\colon P\longmapsto [i]P$ the endomorphism defined
on $E(\F_p)$ which consists of adding $P$ to itself $i$ times. Let  $\pi _{p}\colon E\left(\overline{\F_p}\right) \rightarrow E\left(\overline{\F_p}\right)$, given by $\pi_p(x,y) = \left(x^{p},y^{p}\right)$, be the Frobenius endomorphism on the curve, where $\overline{\F_p}$ is the algebraic closure of the finite field $\F_p$. The relation between the trace $t$ of the Frobenius endomorphism and the group order is given by \cite[Theorem 4.3]{Was08}:
$\sharp E(\F_{p})=p+1-t$ and $\pi_p$ has exactly two eigenvalues $1$ and $p$. This enables us to consider $ P\in \mathbb{G}_{1}=E\left(\overline{\F_p}\right)[r]\cap$ Ker$(\pi _{p}-[1])=E(\F_p)[r]$ and $Q\in \mathbb{G}_{2}=E\left(\overline{\F_p}\right)
[r]\cap$ Ker$(\pi _{p}-[p])$. The ate pairing is defined as follows:
\begin{align*}
e_A\colon \mathbb{G}_{2}\times \mathbb{G}_{1} &\longrightarrow \mu _{r},\\
(Q,P)  &\longmapsto  f_{t-1,Q}(P)^{\frac{p^{k}-1}{r}}.
\end{align*}

In all variants of pairings, one needs a value $f_{m,U}(V)$ which is efficiently computed thanks to the
Miller algorithm \cite{Mil04}. Indeed let $h_{R,S}$  denote a rational function with divisor $ {\rm Div}(h_{R,S})=(R)+(S)-(S+R)-(\mathcal{O})$ where $R$ and $S$ are two arbitrary points on the elliptic curve. In the case of elliptic curves in Weierstrass form, $h_{R,S}=\frac{\ell_{R,S}}{v_{R+S}}$ where $\ell_{R,S}$ is the straight line containing $R$ and $S$ and $v_{R+S}$ is the corresponding vertical line passing through $R+S$.
Miller uses the \textit{double-and-add} method as the addition chains for $m$ (see \cite[Chapter 9]{CohFre06} for more details on addition chains) to compute $f:=f_{m,U}(V)$. Write $m=m_n2^n+...+m_12+m_0>0$ with $m_i\in \{-1,0,1\}$; the (modified) Miller algorithm
that efficiently computes the pairing $f_{m,U}(V)^{(p^k-1)/r}$ of two points $U$ and $V$ is then given as follows:

\medskip

\noindent 1:~~Set $f \leftarrow 1$ and $R \leftarrow U$\\
2:~~\textbf{For} $i = n-1$ \textbf{down to} $0$ \textbf{do}\\
3:\hspace{1.5cm} $f \leftarrow f^{2}\cdot h_{R,R}(V)$, \hspace{1.5cm}  $R \leftarrow 2R$ ~~~~~~~~~~~~~~~~Doubling step\\
5:\hspace{1 cm} \textbf{if} $m_{i}=1$ \mbox{\textbf{then}}\\
6:\hspace{1.5cm} $f \leftarrow f\cdot h_{R,U} (V)$\hspace{1.5cm}$R \leftarrow R + U$, \textbf{end if}~~~~~~Addition step\\
7:\hspace{1 cm} \textbf{if}  $m_{i}=-1$ \mbox{\textbf{then}}\\
8:\hspace{1 cm} $f \leftarrow f/ h_{R,U} (V)$\hspace{1.5cm}$R \leftarrow R - U$,~~\textbf{end for}~~~~~~~Addition step\\
10:~~\textbf{return} $e=f^{\frac{p^{k}-1}{r}}$~~~~~~~~~~~~~~~~~~~~~~~~~~~~~~~~~~~~~~~~~~~~~~~~~~~Final exponentiation

\medskip

\noindent The use of twists enables us to efficiently do some computations during the execution of this algorithm as we explain in the next section.

%%%%%%%%%%%%%%%%%%%%%%
\subsection{Use of Twists}\label{useoftwist}
Twists of elliptic curves enable us to efficiently compute pairings. Indeed, in the Miller algorithm, the doubling of a point (line $3$) and the addition of 
points (lines $6$ and $8$) are done in the extension field $\F_{p^{k}}$ in the case of the ate pairing. The use of twists enables us to perform these operations in a subfield of $\F_{p^{k}}$ and also leads to the denominator elimination. More precisely,
%We observe that the group $\mathbb{G}_{2}$ consists of points with coordinates in $\F_{p^{k}}$. Using the concept of twists of elliptic curve and for efficiency reason in computation, $\mathbb{G}_{2}$ can be represented by a subgroup $\mathbb{G'}_{2}$ of a curve (the twist) isomorphic to $E$ over  $\F_{p^{k}}$.
a twist of an elliptic curve $E$ defined over a finite field $\F_{p}$ is an elliptic curve $E'$ defined over $\F_{p}$ which is isomorphic to $E$ over an algebraic closure of $\F_{p}$. The smallest integer $d$ such that $E$ and $E'$ are isomorphic over $\F_{p^{d}}$ is called the \emph{degree} or the \emph{order }of the twist. Elliptic curves of embedding degrees $k=9$, $15$, $27$ admit twists of order three. Explicit constructions of such curves can be found in \cite{LinZhaZhaWan08}, \cite{DuaCuiCha05} and \cite{BarLynSco02}. The general equation of these curves is given by $E\colon y^{2} = x^{3}+b$.
The equation defining the twist $E'$ has the form $y^{2} = x^{3}+b\omega^{6}$ where $\{1,\omega,\omega^2\}$ is the basis of the $\F_{p^{k/3}}$-vector space $\F_{p^{k}}$ and the isomorphism between $E'$ and $E$ is the map
$\psi\colon E'\longrightarrow E$ given by $\psi(x',y') = (x'/\omega ^{2},y'/\omega^{3})$. Using this isomorphism, points $Q$ in $\mathbb{G}_{2}$ can instead be taken as $(x\omega ^{-2},y\omega^{-3})$ where $(x,y)\in E'(\F_{p^{k/3}})$. The function $h_{R,S}$ is defined by $h_{R,S}(x,y)=\frac{y+\lambda (x_R-x)-y_R}{x-x_{R+S}}$ where $\lambda$ is the slope of the line passing through $R$ and $S$. Observe that using the equation of the curve $y^2=x^3+b$ one has $x-x_{R+S}=\frac{y^2-y_{R+S}^2}{x^2+x_{R+S}x+x_{R+S}^2}$. In the present case of ate pairing, the addition $R+S$ is performed in the extension field $\F_{p^k}$ and the function $h_{R,S}$ is evaluated at a point $(x_P,y_P)\in E(\F_p)$. So, using the twist, the points $R$, $S$ and $R+S$ are taken in the form $(x\omega ^{-2},y\omega^{-3})$ where $(x,y)\in E'(\F_{p^{k/3}})$. Therefore we have $h_{R,S}(x_P,y_P)=\frac{(y_{P} \omega^5+\lambda (x_R \omega^2-x_P \omega^4)-y_R \omega^2)(x_P^2\omega^4+x_{R+S}x_P \omega^2+x_{R+P}^2)}{y_ P^2\omega^9-y_{R+P}^2\omega^{3}}$.\\
We observe that the denominator is an element of the subfield $\F_{p^{k/3}}$ and so will be sent to $1$ during the final exponentiation (line 10 in the Miller algorithm) since $p^{k/3}-1$ is a factor of $(p^k-1)/r$. Consequently we simply ignore that denominator in the Miller algorithm for an efficient computation. More details on twists can be found in \cite{CosLanNae10}.

%%%%%%%%%%%%%%%%%%%%%%%
\subsection{Optimal Pairings}\label{optimalpairing}
The reduction of Miller's loop length is an important way to improve the computation of pairings. The latest work is a generalized method to find the shortest loop, which leads to the concept of optimal pairings due to Vercauteren \cite{Ver10}.
%Indeed, observe that  $r | p^k-1$ but
%$r\nmid p^i-1$ for any $1\leq i< k$. This implies that $r|\Phi_k(q)$ where $\Phi_k$ is the $k-th$ cyclotomic polynomial. Since $T\equiv q$ mod $r$ where $T=t-1$,
%we have $r|\Phi_k(T)$. More generally, if we consider Ate$-i$ pairing, which is a generalisation of Ate
%pairing with Miller function $f_{T_i,Q}$ where $T_i\equiv p^i$ mod $r$, then
%\begin{eqnarray*}r|\Phi_{k/g}(T_i), \mbox{ where } g={\gcd}(i,k)\end{eqnarray*}
%so that the minimal value for $T_i$ is $r^{1/\varphi(k/g)}$ (where $\varphi$ is the Euler's totient function) and the lowest bound is $r^{1/\varphi(k)}$, obtained for $g=1$. The general framework for the computation of optimal pairing is guiven in the following:
%\begin{theo}\cite[Theorem 4]{Ver10}\label{theopt}
%Let $E$ be an elliptic curve defined over $\F_p$. The embedding degree with respect to a large integer $r$ dividing the order of the group $\sharp E(\F_{p})$ is
%denoted $k$.
Let $\lambda=mr$ be a multiple of $r$ such that $r\nmid m$ and write $\lambda=\sum_{i=0}^lc_ip^i=h(p)$, $(h(z)\in \Z[z])$. Recall that $h_{R,S}$ is the Miller function defined in Section \ref{reviewpairing}. For $i=0,\cdots l$ set  $s_i=\sum_{j=i}^lc_jp^j$; then the map
\begin{align}
e_o\colon \mathbb{G}_{2}\times \mathbb{G}_{1} &\longrightarrow \mu _{r} \label{op},\\
(Q,P) &\longmapsto \left(\prod_{i=0}^lf_{c_i,Q}^{p^i}(P)\cdot \prod_{i=0}^{l-1}h_{[s_{i+1}]Q,[c_ip^i]Q}(P)\right)^{\frac{p^k-1}{r}} \nonumber
\end{align}
defines a bilinear pairing and is non-degenerate if $mkp^k \neq ((p^k-1)/r)\cdot \sum_{i=0}^lic_ip^{i-1} \bmod  r$. The coefficients $c_i$, $i=0,\cdots,l$, can be obtained from the short vectors obtained from the lattice
\begin{eqnarray}\label{lattice}
L=\left(\begin{array}{ccccc}
r & 0 & 0 &\cdots &0 \\
-p & 1 & 0 & \cdots &0 \\
-p^2 & 0 & 1 & \cdots &0 \\
\cdots &  \cdots &  \cdots &  \cdots &  \cdots \\
-p^{\phi(k)-1} & 0 & 0 & \cdots &1 \\
\end{array}
\right)
\end{eqnarray}
%\end{theo}

%%%%%%%%%%%%%%%%%%%%%%%%%
\subsection{Final Exponentiation and the Lattice-Based Method for its Computation}\label{final expo and lll}
The result of the Miller loop's step is raised to the power $\frac{p^{k}-1}{r}$. This step is called the final exponentiation (line 10 in Miller's algorithm). The efficient computation of the final exponentiation has became a serious task. Observe that this exponent can be divided into two parts as follows:
$$ \frac{p^{k}-1}{r}=\left[\frac{p^{k}-1}{\phi_k(p)}\right]\cdot\left[\frac{\phi_k(p)}{r}\right]$$
%
%footnote{Note that in certain cases, for efficient computations, $\phi_k(p)$ can be replaced by any order of a cyclotomic subgroup of $\F_{p^k}$: see $k=15.$}$$
where $\phi_k(x)$ is the $k$-th cyclotomic polynomial. The final exponentiation is therefore computed as $f^{\frac{p^{k}-1}{r}}=\left[f^{\frac{p^{k}-1}{\phi_k(p)}}\right]^{\frac{\phi_k(p)}{r}}$. The computation of the first part $A=f^{\frac{p^{k}-1}{\phi_k(p)}}$ is generally inexpensive as it consists of few multiplications, inversion and $p$-th powering in $\F_{p^{k}}$. The second part $A^{\frac{\phi_k(p)}{r}}$ is considered to be more difficult and is called the hard part. An efficient method to compute the hard part is described by Scott \emph{et al.} \cite{ScoBenChaPerKac09a}. They suggested to write $d=\frac{\phi_k(p)}{r}$ in base $p$ as $d=d_0+d_1p+...+d_{\phi (k)-1}p^{\phi(k)-1}$ and find a short vectorial addition chain to compute $A^d$ much more efficiently than the naive method. In \cite{LauKnaRod11}, based on the fact that a fixed power of a pairing is still a pairing, Fuentes \emph{et al.} \cite{LauKnaRod11} suggested to apply Scott \emph{et al.}'s method with a power of any multiple $d'$ of $d$ with $r$ not dividing $d'$. This could lead to a more efficient exponentiation as opposed to computing $A^d$ directly. Their idea for finding the polynomial $d'(x)$ is to apply the $LLL$-algorithm to the matrix formed by $\Q$-linear combinations of the elements $d(x),xd(x),...,x^{\mbox{deg}r-1}d(x)$. They successfully applied this method in the case of elliptic curves of embedding degrees $8$, $12$ and $18$ \cite{LauKnaRod11}. In Sections \ref{k=9} and \ref{k=15} we apply this method to improve the computation of the final exponentiation for elliptic curves of embedding degrees $k=9$ and $15$. A clever method was used by Zhang \emph{et al.} \cite{ZhaLin12} to compute the final exponentiation in the case $k=27$.
%More informations on pairings can be found in \cite{Gal05} and \cite{DuqFre05}.\\

%%%%%%%%%%%%%%%%%%%%%%%%%%%%
\section{Arithmetic in the tower of subfields of $\F_{p^{9}}$, $\F_{p^{27}}$ and $\F_{p^{15}}$}
A pairing is computed as an element of the extension field $\F_{p^k}$. However its efficient computation depends on the arithmetic of subfields of $\F_{p^k}$ which is generally organised as a tower of subfield extensions. In this section we recall the tower extension of finite fields $\F_{p^{9}}, \F_{p^{27}}$ and $\F_{p^{15}}$. We also give explicit cost of the arithmetic operations.
For extension-field arithmetic in $\F_{p^{9}}$ and $\F_{p^{27}}$ we consider $p \equiv 1 \bmod 3$ motivated by the work of Barreto et \emph{al.} \cite{BarLynSco02} on the construction of elliptic curves of embedding degrees $9$ and $27$. This implies that $\F_{p^{k}}$ can be represented as $\F_{p^{k/3}}[X]/(X^3-\alpha)$, for $k=3^i$, $i=1,2,3$, where $\alpha$ is a cubic non-residue modulo $p$.
We choose $p$ such that  $X^3-7$ is irreducible over $\F_p$. Therefore cubic extensions will be constructed using the polynomials $X^3-\alpha_i$ where $\alpha_i=7^{1/3^{i-1}}$. Tower extensions for $\F_{p^{27}}$, together with the one for $\F_{p^{9}}$, are then given by:
\begin{align*}
\F_{p^{3}} &=\F_{p}[u]\text{ with }u^{3}=7 \\
\F_{p^{9}} &= \F_{p^{3}}[v]\text{ with }v^{3}=7^{1/3} \\
\F_{p^{27}} &= \F_{p^{9}}[w]\text{ with }w^{3}=7^{1/9}.
\end{align*}

The costs of the computation of the Frobenius maps and cyclotomic inversions are given in Lemma \ref{frob9}, Lemma \ref{frob27} and Lemma \ref{frob15} for the extensions $\mathbb{F}_{p^{9}}$, $\mathbb{F}_{p^{27}}$, and $\mathbb{F}_{p^{15}}$ respectively. The proof of these lemmas are given in Appendices \ref{arithFP9}, \ref{arithFP27} and \ref{arithFP15}.

\begin{lem}\label{frob9}
	In the finite field $ \mathbb{F}_{p^{9}}$,
	\begin{enumerate}
		\item The computation of the $p^3;p^6$-Frobenius maps costs $6M_{1} $.
		\item The computation of the $p;p^2;p^4;p^5;p^7;p^8$-Frobenius maps costs  $8M_{1}$.
		\item The inverse of an element $\alpha$ of the $G_{\phi_3(p^3)}$-order cyclotomic subgroup is computed as $\alpha^{-1}=\alpha^{p^3}\cdot \alpha^{p^6}$ and the cost is  $36S_1$.
	\end{enumerate}
\end{lem}

Similarly, in the finite field $\F_{p^{27}}$, Lemma \ref{frob27} gives the costs of the computation of the Frobenius maps and cyclotomic inversions. 

\begin{lem}\label{frob27}
	In the finite field $ \mathbb{F}_{p^{27}}$,
	\begin{enumerate}
		\item The computation of the $p^3;p^6;p^9$-Frobenius maps costs $18M_1 $.
		\item The computation of the $p;p^2;p^4;p^5;p^7;p^8$-Frobenius maps costs $26M_1 $.
		\item The inverse of an element $\alpha$ of the $G_{\phi_3(p^9)}$-order cyclotomic subgroup is computed as $\alpha^{-1}=\alpha^{p^9}\cdot \alpha^{p^{18}}$ and the cost is $216S_1$.
	\end{enumerate}
\end{lem}

In the case of $\F_{p^{15}}$, we consider pairing friendly curves over $\F_p$ where $p\equiv 1 \mbox{ mod } 5$ \cite{DuaCuiCha05}. According to
\cite[Theorem 3.75]{Lid97} the polynomial $X^5-\alpha$ is irreducible over $\F_p[X]$ if and only if $\alpha$ is neither a cubic root nor a fifth root in $\F_p$.
A tower extension for $\mathbb{F}_{p^{15}}$ can be constructed as follows:
\begin{align*}
\F_{p^{5}} &=\F_{p}[u]\text{ with }u^{5}=7. \\
\F_{p^{15}} &= \F_{p^{5}}[v]\text{ with }v^{3}=u, \text{ where }u \in \mathbb{F}_{p^{5}}.
\end{align*}
Our main contribution in this section is the computation of Frobenius maps and the inversions in the $\phi_n(.)$-order cyclotomic subgroup of $\F_{p^{k}}^{\ast}$.
The costs of the computation of the Frobenius maps and cyclotomic inversions are given in Lemma \ref{frob15}.

\begin{lem}\label{frob15}
	In the finite field $\F_{p^{15}}$,
	\begin{enumerate}
		\item The computation of the $p^5;p^{10}$-Frobenius maps costs  $10M_1$.
		\item The computation of the $p;p^2;p^3;p^4;p^6;p^7;p^8; p^9$-Frobenius maps costs $ 14M_1$.
		\item The inverse of an element $\alpha$ of the $G_{\phi_3(p^5)}$-order cyclotomic subgroup is computed as $\alpha^{-1}=\alpha^{p^5}\cdot \alpha^{p^{10}}$ and the cost is $54S_1$.
	\end{enumerate}
\end{lem}

In Table \ref{COSTARITHMETICff} we summarise the overall cost of operations in the tower of subfields described above. The costs for squaring, multiplication and inversion are from \cite{LeTan12}, \cite{MraGuiIon09} and \cite{ZhaLin12} respectively for $k=9, 15$ and $27$. Explicit details of the cost of Frobenius maps and inversions in the  cyclotomic subgroups are given in Appendices \ref{arithFP9}, \ref{arithFP27} and \ref{arithFP15}.

\begin{table}
	\centering
	\begin{tabular}{|l|l|l|}
		
		\hline
		Fields & Operations& Costs\\
		
		\hline
		\multirow{3}{*}{$\F_{p^3}$}
		%\cline{2-3}
		&Multiplication $M_3$ & $6M_1$\\
		&Squaring $S_3$& $5S_1 $ \\
		&Inversion $I_3$& $I_1 + 9M_1 +2S_1$\\
		\hline
		\multirow{3}{*}{$\F_{p^9}$}
		&Multiplication $M_9$& $36M_1$\\
		&Squaring $S_9$&  $25S_1$\\
		&Inversion $I_9$ &  $ I_1+ 63M_1 +12S_1$\\
		&Frobenius $p^3;p^6$& $ 6M_1$\\
		&Frobenius $p;p^2;p^4;p^5;p^7;p^8$&  $8M_1$\\
		& Inversion in $G_{\phi_3(p^3)}$& $18M_1+15S_1$\\

		\hline
		\multirow{3}{*}{ $\F_{p^{27}}$}
		
		&Multiplication $M_{27}$& $216M_1$\\
		&Squaring $S_{27}$&$ 125S_1$\\
		&Inversion $I_{27}$ & $ I_1 + 387M_1+62S_1$\\
		&Frobenius $p^3;p^6;p^9$& $18M_1 $\\
		&Frobenius $ p;p^2;p^4;p^5;p^7;p^8$ & $26M_1$\\
		& Inversion in $G_{\phi_3(p^9)}$& $108M_1+75S_1$\\
		\hline
		\multirow{3}{*}{ $\F_{p^5}$}
		& Multiplication $M_5$ &$ 9M_1$\\
		&Squaring $S_5$&$ 9S_1$\\
		&Inversion $I_5$ & $ 1I_1+45M_1 + 5S_1$\\
		\hline
		\multirow{3}{*}{ $\F_{p^{15}}$}
		&Multiplication $M_{15}$& $45M_1$\\
		&Squaring $S_{15}$& $45S_1$ \\
		&Inversion $I_{15}$ & $I_1+126M_1 + 23S_1$ \\
		&Frobenius $p^5;p^{10}$& $10M_1$ \\
		&Frobenius $ p;p^2;p^3;p^4;p^6;p^7;p^8; p^9$ &$ 14M_1$\\
		& Inversion in $G_{\phi_3(p^5)}$& $27M_1+27S_1$\\
		\hline
	\end{tabular}
	\caption{Cost of operations in extension fields from \cite{LeTan12}, \cite{MraGuiIon09} and \cite{ZhaLin12} and this work (see Appendices \ref{arithFP9}, \ref{arithFP27} and \ref{arithFP15})}\label{COSTARITHMETICff}
\end{table}

%\begin{rem}	
%In the next sections, the notations $m_c, s_c$ will be used. They denote the costs of a multiplication and a squaring in the field $\F_p$ where the bit size of $p$ is $c$ respectively. These notations are not really different from $M_1$ and $S_1$ but here we insist on the size of the base field $\F_{p}$.
%\end{rem}

%%%%%%%%%%%%%%%%%%%%%%%%%%
\section{Elliptic Curves with Embedding Degree 9}\label{k=9}

This section describes the computation of the optimal ate pairing (Miller step and the final exponentiation) on the parameterized elliptic curve defined in \cite{LinZhaZhaWan08}. The correctness of the results can be verified with the MAGMA code available in \cite{FouNadAmi18}. This family of elliptic curves has embedding degree 9 and a $\rho$-value $1.33$ and is parameterized by :
\begin{align*}
p &= ((x+1)^2+((x-1)^2(2x^3+1)^2)/3)/4,\\
r &=(x^6+x^3+1)/3,\\
t &=x+1.
\end{align*}

%%%%%%%%%%%%%%%%%%%%%%
\subsection{Optimal ate pairing}

Based on the general framework described by Vercauteren in \cite{Ver10}, the short vector obtained from the lattice
$L$ defined by equation (\ref{lattice}) gives the optimal function $h(z)=\sum_{i=0}^5c_iz^i=x-z\in \Z[z].$ A straightforward application of formula (\ref{op}) yields the optimal pairing
\begin{align*}
e_o\colon \mathbb{G}_{2}\times \mathbb{G}_{1} &\longrightarrow \mu_{r} ,\\
(Q,P) &\longmapsto f_{x,Q}(P)^{\frac{p^9-1}{r}}.
\end{align*}

%%%%%%%%%%%%%%%%%%%
\subsection{Cost of the execution of the Miller loop}

The Miller loop consists of the doubling steps (line 3 in the Miller algorithm) and addition steps (line 6 or 8 in the Miller algorithm). These steps use the Miller function $h_{R,S}$ either in affine coordinates or in projective coordinates. The work of Zhang \emph{et al.} \cite[Section 3]{ZhaLin12} presents the currently fastest formulas in projective coordinates.  The doubling step costs $9M_1+3M_3+9S_3$ and the cost of the addition step is $9M_1+12M_3+5S_3$. For an explicit cost of the computation of $f_{x,Q}(P)$, we wrote a Pari/GP code to find a suitable $x$ with low Hamming weight and minimal number of bits for the 128 bit-security level according to Table \ref{recomendationParameters}. The best value we were able to find is $x=2^{43}+2^{37}+2^{7}+1$ which gives an $r(x)$ prime of 257 bits and $p(x)$ a prime of 343 bits. The values $p$ and $x$ are both congruent to 1 modulo 6 so that the corresponding elliptic curve is $y^2=x^3+1$ \cite{LeTan12}. The computation of $f_{x,Q}(P)$ therefore costs 43 doubling steps, 3 additions, 42 squaring and 45  multiplications in $\F_{p^{9}}$. Thus the total cost for the computation of the Miller loop for the optimal pairing on elliptic curves of embedding degree 9 is $43(9M_1+3M_3+9S_3)+3(9M_1+12M_3+5S_3)+42S_9+45M_9$, that is, $45M_9+165M_3+414M_1+42S_9+402S_3$. Using the arithmetic in Table \ref{COSTARITHMETICff}, the overall cost is $3024M_1+3060S_1$. As far as we are aware, no other explicit cost with a specific value of $x$ is reported in the literature.

%%%%%%%%%%%%%%%%%%%%%%%
\subsection{Cost of the computation of the final exponentiation}

As explained in Section \ref{intro}, the final exponentiation in this case can be divided as
$$f^{(p^9-1)/r}=\left(f^{p^3-1}\right)^{ (p^6+p^3+1)/r}=\left(f^{p^3-1}\right)^ d.$$
 We then used the lattice method described by Fuentes \emph{et al.} \cite{LauKnaRod11} that we briefly explained in Section \ref{final expo and lll}. It is applied to the following matrix in which the coefficient 243 is used to obtain integer entries as $d=(p^6+p^3+1)/r$ is a polynomial with rational coefficients  \begin{eqnarray}
M=\left(
\begin{array}{r}
243d(x) \\
243xd(x) \\
243x^2d(x) \\
243x^3d(x)  \\

243x^4d(x)  \\
243x^5d(x)  \\
\end{array}
\right)\end{eqnarray}
We obtain the following multiple of $d$: $d'=x^3d=k_0+k_1p+k_2p^2+k_3p^3+k_4p^4+k_5p^5$ where the polynomials $k_i$, $i=0,\dots,5$ are as follows
\begin{eqnarray*}
	\begin{array}{lll}
		k_0=-x^4 + 2x^3 - x^2,~~~& k_1=-x^3 + 2x^2 - x,~~~& k_2=-x^2 + 2x - 1,\\
		k_3=x^7 - 2x^6 + x^5 + 3,~~~&k_4= x^6 - 2x^5 + x^4, ~~~& k_5= x^5 - 2x^4 + x^3.
	\end{array}
\end{eqnarray*}
They verify the relations (see the code in \cite{FouNadAmi18} for verification)
$$k_2=-(x-1)^2,\quad k_1=xk_2, \quad k_0=xk_1, \quad k_5=-xk_0, \quad k_4=xk_5, \quad k_3=xk_4+3.$$
If we set $A=f^{p^3-1}$ then
\begin{itemize}
	\item The cost for the computation of $A$ is one $p^3$-Frobenius, one inversion in $\F_{p^{9}}$ and one multiplication in $\F_{p^9}$.
	\item The cost of the computation of $A^{k_0}$, $A^{k_1}$ and $A^{k_4}$ is three exponentiations by $x$,
	\item The cost of the computation of $A^{k_5}$ is one inversion in the cylotomic subgroup and one exponentiation by $x$.
	\item The cost of the computation of $A^{k_2}$ is one inversion in the cyclotomic subgroup and two exponentiations by $(x-1)$.
	\item The cost of the computation of $A^{k_3}$ is two multiplications, one squaring and one exponentiation by $x$.
	%\item The product $A^{k_0}$ \cdot$A^{k_1}$ \cdot$A^{k_2}$ \cdot$A^{k_3}$ \cdot$A^{k_4}$ \cdot$A^{k_5}$ consists of $5$ multiplications.
\end{itemize}

Note that the inversion in the cyclotomic subgroup $G_{\phi_3(p^3)}$ of order $p^6+p^3+1$ is computed as $A^{-1}=A^{p^3}\cdot A^{p^6}$ (see Appendix \ref{arithFP9} for details and cost). The cost for the hard part $A^{d'}$ is then two exponentiations by $x-1$, five exponentiations by $x$, seven multiplications in $\F_{p^9}$, one squaring in $\F_{p^{9}}$, two cyclotomic inversions $I_{G_{\phi_3(p^3)}}$ and $p,p^2, p^3,p^4,p^5$-Frobenius maps. Using the value of $x$ given above, one exponentiation by $x$ costs $43S_9+3M_9$ whereas one exponentiation by $x-1$ costs $43S_9+2M_9$. Finally the hard part costs
$$2(43S_9+2M_9)+5(43S_9+3M_9)+7M_9+1S_9+2I_{G_{\phi_3(p^3)}}=302S_9+26M_9+2I_{G_{\phi_3(p^3)}}$$
and $p,p^2, p^3,p^4,p^5$-Frobenius maps. The total cost of the final exponentiation is $ 1I_9+27M_9+302S_9+2I_{G_{\phi_3(p^3)}}$ and $p,p^2, 2*p^3,p^4,p^5$-Frobenius maps.

%%%%%%%%%%%%%%%%%%%%%%%
\subsection{Improvement and comparison with previous work}

From the results in \cite{LeTan12}, the hard part costs $309S_9+50M_9$ and $p,p^2, p^3,p^4,p^5$-Frobenius maps. If we include the cost $ 1I_9+1M_9$ and $p^3$-Frobenius for the easy part and using the arithmetic in Table \ref{COSTARITHMETICff}, the overall cost is $I_1+1115M_1+7592S_1$ for this work as opposed to $I_1+1943M_{1}+7737S_1$ for Le et \emph{al.} \cite{LeTan12}. We therefore save $828M_1+145S_1$ compared to their work. Although the sizes of $p$ are $343$ bits in this work and $348$ bits in \cite{LeTan12}, the cost of a  multiplication is almost the same in the two corresponding fields (see Section \ref{compare}).

%%%%%%%%%%%%%%%%%%%%%%%%
\section{Elliptic Curves with Embedding Degree 15}\label{k=15}

In this section we give explicit formulas; together with their cost, for the Miller loop in the computation of the optimal ate pairing. We then compute the cost of the final exponentiation on the parameterized elliptic curve defined in \cite{DuaCuiCha05}. The correctness of the results can be verified in \cite{FouNadAmi18}. This family of elliptic curves has embedding degree 15 and a $\rho$-value $1.5$ and is parameterized by :
\begin{align*}
p &=(x^{12}-2x^{11}+x^{10}+x^7-2x^6+x^5+x^2+x+1)/3,\\
r &=x^8-x^7+x^5-x^4+x^3-x+1,\\
t &=x+1.
\end{align*}

%%%%%%%%%%%%%%%%%%%%%%%%%%%%%%
\subsection{Optimal ate pairing}

The Vercauteren approach described in \cite{Ver10} enabled us to obtain the short vector  from the lattice
$L$ defined by Equation (\ref{lattice}) which lead to the optimal function $h(z)=\sum_{i=0}^5c_iz^i=x-z\in \Z[z].$ A straightforward application of formula (\ref{op}) yields the optimal pairing
\begin{align*}
e_o\colon \mathbb{G}_{2}\times \mathbb{G}_{1} & \longrightarrow \mu_{r},\\
(Q,P) &\longmapsto f_{x,Q}(P)^{\frac{p^{15}-1}{r}}.
\end{align*}

%%%%%%%%%%%%%%%%%%%%%%%
\subsection{Cost of the computation of the Miller loop} 

In this section, we consider the Miller function given in affine coordinates, following the analysis of Lauter et \emph{al.} \cite{LauMonNae10} who suggested the use of affine coordinates at a higher security level. The Miller function used for the computation of $f_{x,Q}(P)$ in this case is described in \cite[Table 2]{ZhaLin12} with the fastest cost to date. At 192-bits security level on elliptic curves with $k=15$, the best value of $x$ we were able to find with a Pari/GP code is $x=2^{48}+2^{41}+2^{9}+2^8+1$. This value gives an $r(x)$ prime of $385$ bits and $p(x)$ of 575 bits which correspond to parameters for 192-bit security level according to Table \ref{recomendationParameters}. The value of $p$ is congruent to 1 modulo 5 and a curve equation can be $y^2=x^3+1$. The Miller loop consists here of computing $f_{x,Q}$ which costs 48 doubling steps, 4 additions steps, 47 squaring and 51  multiplications in $\F_{p^{15}}$. Considering the currently fastest cost for doubling and addition step in \cite[Table 2]{ZhaLin12}, the Miller loop costs $48(15M_1+3M_5+2S_5+I_5)+4(15M_1+3M_5+2S_5+I_5)+47S_{15}+51M_{15}$, that is, $51M_{15}+156M_5+780M_1+47S_{15}+104S_5+52I_5$. Using the arithmetic in Table \ref{COSTARITHMETICff}, the overall cost is $52I_{1}+6819M_{1}+3311S_{1}$. As far as we are aware, no explicit cost is reported in the literature in the case $k=15$ with a specific value of $x$.

%%%%%%%%%%%%%%%%%%%%%%%%%
\subsection{Cost of the computation of the final exponentiation}

The final exponentiation in this case is written in a different way as $f^{(p^{15}-1)/r}=\left(f^{p^5-1}\right)^{ (p^{10}+p^5+1)/r}=\left(f^{p^{5}-1}\right)^ d$. This decomposition is used instead of $ \frac{p^{15}-1}{r}=\left[\frac{p^{15}-1}{\phi_{15}(p)}\right]\cdot\left[\frac{\phi_{15}(p)}{r}\right]$ as usually done, for efficiency reasons in the computation. Observe that $\frac{p^{15}-1}{\phi_{15}(p)}=p^7+p^6+p^5-p^2-p-1$ and $\phi_{15}(p)=p^8-p^7+p^5-p^4+p^3-p+1$ will lead to several multiplications and Frobenius map operations. Thus the lattice method described by Fuentes \emph{et al.} \cite{LauKnaRod11} that we briefly explained in Section \ref{final expo and lll} is applied to the following matrix $M$.  In the matrix $M$ the coefficient $\lambda = \frac{59049}{19683}$ is used to obtain integer entries as $d=(p^{10}+p^5+1)/r$ is a polynomial with rational coefficients.
\begin{eqnarray*}
M=\left(
\begin{array}{c}
\lambda d(x) \\
\lambda xd(x) \\
\lambda x^2d(x) \\
.  \\
.\\
\lambda x^7d(x)  \\
\end{array}
\right)\end{eqnarray*}
We then obtained the following multiple of $d$: $d'=3x^3d=k_0+k_1p+...k_9p^9$, where the polynomials $k_i, i=0,...,9$ are defined as follows

\begin{eqnarray*}
	\begin{array}{lcl}
		k_0=-x^6 + x^5+x^3- x^2, &\quad k_1=-x^5 + x^4+x^2 - x,&\\
		k_2=-x^4 +x^3+ x - 1,&&\\
		\multicolumn{2}{l}{k_3=x^{11} - 2x^{10} + x^9 + x^6 - 2x^5 + x^4 - x^3 + x^2 + x + 2,}&\\
		k_5=x^{11} - x^{10} - x^8 + x^7 + 3,&&\\
		\multicolumn{2}{l}{k_4=x^{11} - x^{10} - x^9 + x^8 + x^6 - x^5 - x^4 + x^3 - x^2 + 2x + 2,}&\\
		k_6=x^{10} - x^9 - x^7 + x^6, &\quad k_7=x^9 - x^8 - x^6 + x^5,&\\
		k_8=x^8 - x^7 - x^5 + x^4,& \quad k_9=x^7 - x^6 - x^4 + x^3.&\\
	\end{array}
\end{eqnarray*}
The polynomials $k_i$, $i=0,\dots, 9$ verify the relations (see the code in \cite{FouNadAmi18} for verification)
\begin{eqnarray*}
	\begin{array}{l}
		k_2=-(x-1)^2(x^2+x+1),\quad k_1=xk_2, \quad k_0=xk_1,\\
		k_9=-xk_0, \quad k_8=xk_9,\quad k_7=xk_8,\\
		k_6=xk_7, \quad k_5=xk_6+3,\quad k_4=M-(k_1+k_7),\\
		k_3=M-(k_0+k_6+k_9) \mbox{ where } M=(k_2+k_5+k_8).
	\end{array}
\end{eqnarray*}
Set $A=f^{p^5-1}$; then
\begin{itemize}
	\item The cost for the computation of $A$ is one $p^5$-Frobenius, one  inversion in $\F_{p^{15}}$ and one multiplication in $\F_{p^{15}}$.
	\item The computation of $A^{k_2}$ is two exponentiations by $x$, two exponentiations by $x-1$, two multiplications and one cyclotomic inversion.
	\item The cost of the computation of  $A^{k_0}$,$A^{k_1}$,$A^{k_6}$,$A^{k_7}$ is five exponentiations by $x$; the computation of $A^{k_9}$ costs one exponentiation by $x$ and one cyclotomic inversion.
	\item The computation of $A^{k_5}$ is one exponentiation by $x$, two multiplications and one squaring in $\F_{p^{15}}$.
	\item The computation of $A^{k_4}$ costs four multiplications in $\F_{p^{15}}$ and one cyclotomic inversion.
	\item The computation of $A^{k_3}$ costs three multiplications in $\F_{p^{15}}$ and one cyclotomic inversion.
	%\item The product $A^{k_0}$ \cdot$A^{k_1}$ \cdot$A^{k_2}$ \cdot$A^{k_3}$ \cdot$A^{k_4}$ \cdot$A^{k_5}$
\end{itemize}
Therefore, the cost of the computation of $A^{d'}$ is two exponentiations by $x-1$, nine exponentiations by $x$, 20 multiplications, one squaring in $\F_{p^{15}}$, four inversions in the cyclotomic subgroup $G_{\phi_3(p^5)}$ of order $p^{10}+p^5+1$ (note that $A^{-1}=A^{p^5}\cdot A^{p^{10}}$; see Appendix \ref{arithFP15} for details) and $p,p^2,p^3,p^4,p^5,p^6,p^7,p^8,p^9$-Frobenius maps.  Using the value of $x$ given above, the cost of the hard part is $2(48S_{15}+3M_{15})+9(48S_{15}+4M_{15})+20M_{15}+1S_{15}+4I_{G_{\phi_3(p^5)}}=529S_{15}+62M_{15}+4I_{G_{\phi_3(p^5)}}$ and $p,p^2,p^3,p^4,p^5,p^6,p^7,p^8,p^9$-Frobenius maps.
% The cost for the computation of $A$ is $1$ Inversion in $\F_{p^{15}}$, $1$ multiplication in $\F_{p^{15}}$ and $1$ $p^5$-Frobenius. The cost for the computation of $A$ is $1$ Inversion in $\F_{p^{15}}$, $1$ multiplication in $\F_{p^{15}}$ and $1$ $p^5$-Frobenius.
The total cost of the final exponentiation in this work is therefore $1I_{15}+529S_{15}+63M_{15}+4I_{G_{\phi_3(p^5)}}$ and $p,p^2,p^3,p^4,2*p^5,p^6,p^7,p^8,p^9$-Frobenius maps.

\begin{rem}\label{rem}
The cost given by Le \emph{et al.} \cite{LeTan12} for the hard part is $11$ exponentiations by $x$, 22 multiplications, 2 inversions in $\F_{p^{15}}$ and $9$ Frobenius maps. The authors said that the cost of an inversion in $\F_{p^{15}}$ is free; with a reference to a similar computation, but on elliptic curves with even embedding degree, unfortunately we do not see how this is possible. Also, they considered an $x$ of 64 bits and Hamming weight 7 and claimed that the cost is $88M_{15}+528S_{15}$ instead of $11(6M_{15}+64S_{15})=88M_{15}+704S_{15}$. Therefore, if we count the 2 inversions in $\F_{p^{15}}$ (these inverses are in fact in the cyclotomic subgroup $G_{\phi_3(p^5)}$), then their final cost is $88M_{15}+704S_{15}+2I_{G_{\phi_3(p^5)}}$ and $11$ Frobenius maps, whereas our cost is $ 62M_{15} +529S_{15}+4I_{G_{\phi_3(p^5)}}$.
\end{rem}

%%%%%%%%%%%%%%%%%%%%%%%%%%
\subsection{Improvement and comparison with previous work}

Considering Remark \ref{rem}, the cost of the final exponentiation in \cite{LeTan12} is $1I_{15}+704S_{15}+89M_{15}+2I_{G_{\phi_3(p^5)}}$ and $p,p^2,p^3,p^4,2*p^5,p^6,p^7,p^8,p^9$-Frobenius maps. We observe that we have improved the results by saving $26M_{15}+175S_{15}-2I_{G_{\phi_3(p^5)}}$.
Using the arithmetic in Table \ref{COSTARITHMETICff}, the overall cost is $I_1+3066M_{1}+24071S_1$ for this work and $I_1+4236M_1+31838S_1$ for Le et \emph{al.} \cite{LeTan12}. We therefore save $26M_{15}+175S_{15}-2I_{G_{\phi_3(p^5)}}=1170M_1+7767S_1$ as compared to their work. A MAGMA code for the implementation to ensure the correctness of the decomposition of the final exponentiation and the Miller function is available in \cite{FouNadAmi18}.

%%%%%%%%%%%%%%%%%%%%%%%%%
\section{Elliptic Curves with Embedding Degree 27}\label{k=27}

The parameterized elliptic curve with embedding degree $27$ is defined in \cite{BarLynSco02}. This family has  a $\rho$-value $10/9$ and is parameterized by the following polynomials:
\begin{align*}
p &=1/3(x-1)^2(x^{18}+x^9+1)+x,\\
r &=1/3(x^{18}+x^9+1),\\
t &=x+1.
\end{align*}

%%%%%%%%%%%%%%%%%%%%
\subsection{The Miller loop and the final exponentiation}

The Miller loop  and the final exponentiation were  studied by Zhang and Lin in \cite{ZhaLin12}. They found the optimal function
$h(z)=\sum_{i=0}^{17}c_iz^i=x-z\in \Z[z]$ and the optimal pairing is given by
\begin{align*}
e_o\colon \mathbb{G}_{2}\times \mathbb{G}_{1} &\longrightarrow \mu_{r} ,\\
(Q,P) &\longmapsto f_{x,Q}(P)^{\frac{p^{27}-1}{r}}.
\end{align*}
Zhang and Lin used the parameter $x=2^{28}+2^{27}+2^{25}+2^{8}-2^3$ for their computation at $256$-bit security level.
The cost of the Miller step that they obtained is therefore $28(3M_9+2S_9+1I_9+9M_1)+4(3M_9+2S_9+1I_9+9M_1)+27(6S_{9})+30(6M_{9})+1I_{27}=276 M_9+226S_9+32I_9+288M_1+I_{27}$ operations. The computation of the final exponentiation in \cite{ZhaLin12} requires $1I_{27}+11M_{27}$, $17$ powers of $x$, $2$ powers of $x-1$ and $p,p^2,p^3,p^4,p^5,p^6,p^7,p^8,2*p^9$-Frobenius maps. Therefore the explicit cost of the final exponentiation is  $1I_{27}+17(4(6M_9)+28(6S_9)+36M_1)+2(5(6M_9)+28(6S_9)+36M_1)+11(6M_{9})+228M_{1}=1I_{27}+648 M_9+3192S_9+912M_1$.

Then the explicit cost for the computation of the Miller loop and the final exponentiation given in that work is $12627M_1+8670S_1+33I_1$ and $24627M_1+114998S_1+1I_1$ respectively (see   \cite{ZhaLin12} for details).

\begin{rem}
	The negative coefficient in the value of $x$ affects the efficiency since one full inversion in $\F_{p^{27}}$ is required in the Miller algorithm (line 8)  and also $19$ inversions in the cyclotomic subgroup are required when raising to the power of $x$ during the final exponentiation.
\end{rem}

In the next section we explain the choice of another parameter to avoid these additional operations.
%\begin{rem}
%We observed that the authors in \cite{ZhaLin12} said that the inversion of $M_0\in\G_3(p^9)$ when exponetiating to $x$ is computed as $M_0^{p^{18}}+M_0^{p^9}$ but we think that this is a multiplication. Also the computation of the easy part requires a multiplication. Therefore their cost will be added of $20$ more multiplications in $\F_{p^{27}}$.
%\end{rem}

%%%%%%%%%%%%%%%%%%%%%%%%%
\subsection{Improvement and comparison with previous work}

We use the arithmetic (especially  the computation of inversion in the cyclotomic subgroup) and a specific value of $x$ to improve the costs in \cite{ZhaLin12}. To be more precise, a careful search with a Pari/GP code enabled us to find the value $x=2^{29}+2^{19}+2^{17}+2^{14}$ so that $r$ has a prime factor of length $514$ bits and the prime $p$ has a bit length of 579 for 256-bit security level according to Table \ref{recomendationParameters}. An adequate elliptic curve has the equation $y^2=x^3-2$. The corresponding base field is of length 579 bits which is a little larger than 573, the length of the base field in \cite{ZhaLin12}. However $m_{579}\approx m_{573}$, (see section \ref{compare} for notations) and we have an extra doubling step; we avoid the full inversion in $\F_{p^{27}}$ and $17$ inversions in the cyclotomic subgroup $\G_{\phi_3(p^9)}$ when raising to power $x$. We perform two inversions in the cyclotomic subgroup only when raising to power $x-1$. The cost of the Miller loop now becomes $29(3M_9+2S_9+1I_9+9M_1)+3(3M_9+2S_9+1I_9+9M_1)+27(6S_{9})+30(6M_{9})=276M_9+226S_9+32I_9+288M_{1}$. Using the arithmetic in Table \ref{COSTARITHMETICff}, the overall cost for the Miller loop is $32I_1+12240M_1+6034S_1$ for this work where we saved at least one inversion in $\F_{p^{27}}$.

Our cost for the final exponentiation is $1I_{27}+17(3(6M_9)+29(6S_9))+
2(4(6M_9)+29(6S_9))+2I_{\G_{\phi_3(p^9)}}+11(6M_{9})=1I_{27}+420 M_9+3306S_9+2I_{\G_{\phi_3(p^9)}}$  and $p,p^2,p^3,p^4,p^5,p^6,p^7,p^8,2*p^9$-Frobenius maps. Using the arithmetic in Table \ref{COSTARITHMETICff}, the overall cost is $I_1+15951M_1+82862S_1$ for this work. The implementation of this pairing is available in \cite{FouNadAmi18}.

%\begin{rem}\label{compararith}
%Observe that we are performing arithmetic over a prime fields of size $579$ bits while the work of Zhang \emph{et al.}\cite{ZhaLin12} is based on a prime field of size $573$ bits. In view of giving a clear comparaison with our work we compare the arithmetic in the two fields. Asuming that arithmetic is done ion a $64$- bits platform, Aranha \emph{et al.} \cite[Section 8]{AraFueMen12} mentioned that a finite field element in characteristic $p$ can be represented with $m=1+\mbox{log}_2(p)$ binary coefficients stored in $n_{64}=\lceil\frac{m}{64}\rceil $ $64$-bit processor word, and a multiplication is done with a complexity $\mathcal{O}(2n_{64}^2+n_{64}$. Thi senables us to make the approximation $1M_1=1.25M_1$.
%\end{rem}
%Cosidering the remark \ref{compararith}, the cost of the Miller step for this work is

%%%%%%%%%%%%%%%%%%%%%%%%
\section{General Comparison  }\label{fcomparaison}

In this section, we summarize the different costs obtained in this work and compare our results with previous works.

\begin{table}[h]
	\centering
	\footnotesize
	\begin{tabular}{|c|l|c|c|}
		
		\hline
		\textbf{Curves} & \textbf{References} & \textbf{Miller loop} &\textbf{Final Exponentiation} \\
		
		\hline
		\multirow{2}{*}{ \large $k=9 \atop 128\text{-bit}$}&Previous work \cite{LeTan12} & \textit{No specific cost reported} &$I_1+1943M_1+7737S_1$\\
		&This work&$3024M_1+3060S_1$&$I_1+1115M_1+7592S_1$\\
		
		\hline		
		\multirow{2}{*}{\large $k=15 \atop 192\text{-bit}$}&Previous work \cite{LeTan12} & \textit{No specific cost reported}  &$I_1+4263M_1+31811S_1$\\
		& This work & $52I_1+6819M_1+3311S_1$     &$I_1+3093M_1+24044S_1$\\
		
		\hline
		
		\multirow{2}{*}{\large $k=27 \atop 256\text{-bit}$}  &  Previous work \cite{ZhaLin12} & $33I_1+ 12627M_1+8670S_1$& $I_1+24627M_1+114998S_1$\\
		& This work & $32I_1+12240M_1+6034S_1$    &$I_1+15951M_1+82862S_1$\\
		\hline
	\end{tabular}
	\caption{Comparison of the cost of the  Miller loop and the final exponentiation.}\label{COSTCOMARISON}
\end{table}

If we assume that the cost of a squaring is the same as the cost of a multiplication then the cost of the final exponentiation is $ I_1+27137M_1$ and $ I_1+36074M_1$ for this work and previous work \cite{LeTan12} respectively for $k=15$. The theoretical improvement obtained in this work is therefore up to 25\%. A similar analysis with $k=9$  yields an improvement of 8\%. In the case of  curves with $k=27$, our results present an  improvement of 14\% and 29\% for the Miller loop and final exponentiation respectively as compared to the work in \cite{ZhaLin12}.

%%%%%%%%%%%%%%%%%%%%%%%%
\section{New Parameters for Optimal Ate pairing on Elliptic Curves with embedding degrees 9, 15 and 27}\label{update}

In this section we consider new parameters for parameterized curves of embedding degrees 9, 15 and 27 at the 128-bit, 192-bit and 256-bit security levels. We consider recent advances in the computation of the discrete logarithm thanks to the Number Field Sieve (NFS) algorithm and its variants described in some papers as mentioned in the introduction. Mostly, the paper of Barbulescu and Duquesne \cite{DuqBar18} presents a detailed approach for computing new parameters for pairings. Scott and Guillevic \cite{ScoGui18} have proposed tentative general recommended parameters for classical security level and that we reproduce here in Table \ref{RNP}.

\begin{table}[h]
	\centering
	\begin{tabular}{|c|c|c|c|c|}
		\hline
		% after \\: \hline or \cline{col1-col2} \cline{col3-col4} .
		DL Algorithm&  AES-$128$ & AES-$192$   & AES-$256$  \\
		%	Algorithm  &    $r$       &    $p^k$            &    $\rho\thickapprox 1$ &    $\rho\thickapprox 2$   \\
		%&$p^k$ & Bit length of\\$r$ & $k$, $\rho\thickapprox 1$& $k$, $\rho\thickapprox 2$\\
		\hline
		
		NFS & $3072$ & $7680$ & $15360$ \\
		\hline
		exTNFS& $3618$ & $9241$ & $18480$\\
		\hline
		SexTNFS & $5004$& $12871$ & $27410$\\
		\hline
	\end{tabular}
	\caption{Recommended extension fields size ($p^k$) to obtain desired levels of security \cite{ScoGui18}.}\label{RNP}
\end{table}

Following Table \ref{RNP}, we searched for new parameters that will ensure resistance to SexTNFS algorithm at the various security levels for curves of embedding degrees 9, 15 and 27.

%%%%%%%%%%%%%%%%%%%%%%%%%%
\subsection{New parameters and costs for optimal ate pairing at the 128-bit security level for $k=9$ and $k=15$}

\begin{itemize}
	\item Case of $k=9$. Following the recommendation from Table \ref{RNP}, we found the value $x=2^{70}+2^{59}+2^{46}+2^{41}+1$. This gives a prime $p$ of 559 bits and a prime $r$ of 419 bits. We proceed as described in Section \ref{k=9} to obtain the cost of the Miller loop and the final exponentiation. The Miller loop in this case costs $70(9M_1+3M_3+9S_3)+4(9M_1+12M_3+5S_3)+69S_9+73M_9$. This is equal to $ 73M_9+258M_3+666M_1+69S_9+650S_3$. Using the arithmetic in Table \ref{COSTARITHMETICff}, the overall cost is $4842M_1+4975S_1$.  Using the value of $x$ given above,  the hard part of the final exponentiation costs  $2(70S_9+3M_9)+5(70S_9+4M_9)+7M_9+1S_9+2I_{G_{\phi_3(p^3)}}=491S_9+33M_9+2I_{G_{\phi_3(p^3)}}$ and $p,p^2, p^3,p^4,p^5$-Frobenius maps. The total cost of the final exponentiation is $ 1I_9+34M_9+491S_9+2I_{G_{\phi_3(p^3)}}$ and $p,p^2, 2*p^3,p^4,p^5$-Frobenius maps for a total cost of $ I_1+1367M_1+12317S_1$.
	\item Case of $k=15$. Following the recommendation from Table \ref{RNP}, we found the value $x=2^{31}+2^{19}+2^{5}+2^{2}$. This gives a prime $p$ of 371 bits and a prime $r$ of 249 bits which is close to 256 bits as required to have 128 bits on the curve side. We proceed as described in Section \ref{k=15} to obtain the cost of the Miller loop and the final exponentiation. The Miller loop in this case costs $3(15M_1+13M_5+3S_5)+31(15M_1+6M_5+7S_5)+30S_{15}+33M_{15}$. This is equal to $ 33M_{15}+225M_5+510M_1+30S_{15}+226S_5$. Using the arithmetic in Table \ref{COSTARITHMETICff}, the overall cost is $4020M_1+3384S_1$.  Using the value of $x$ given above,  the hard part of the final exponentiation costs  $2(31S_{15}+4M_{15}+1I_{G_{\phi_3(p^5)}})+9(31S_{15}+3M_{15})+20M_{15}+1S_{15}+4I_{G_{\phi_3(p^5)}}=55M_{15}+342S_{15}+6I_{G_{\phi_3(p^5)}}$ and $p,p^2, p^3,p^4,p^5,p^6, p^7, p^8, p^9$-Frobenius maps. The total cost of the final exponentiation is $ I_{15}+56M_{15}+342S_{15}+6I_{G_{\phi_3(p^5)}}$ and $p,p^2, p^3,p^4,2*p^5,p^6, p^7, p^8, p^9$-Frobenius maps for a total cost of $ I_1+2940M_1+15575S_1$.
\end{itemize}

Table \ref{compare128} below compares our results with previous results at the 128-bit security level.
\begin{center}
	\footnotesize
	\begin{table}[h]
		\begin{tabular}{|c|c|c|c|c|c|}
			
			\hline
			\textbf{Curves-Ref.} & \textbf{Miller loop} &\textbf{Final Exp.} &Size of $p$&Total($S_1=M_1$)\\
			
			\hline
			$k=15$ & $4020M_1+3384S_1$ &$ I_1+2940M_1$&$371$&$I_1+25919M_1$\\
			(This work) &&$+15575S_1$&&\\
			\hline
			KSS$16$ \cite{DuqBar18}&$7534M_1$&$I_1+18542M_1$&$340$&$I_1+26076M_1$\\
			\hline
			BLS$12$ \cite{DuqBar18}&$7708M_1$&$I_1+8295M_1$&$461$&$I_1+16003M_1$\\
			\hline
			BN$12$ \cite{DuqBar18,AraFueMen12}&$12068M_1$&$I_1+7485M_1$&$461$&$I_1+19553M_1$\\
			\hline
			$k=9$&$4842M_1+4975S_1$&$ I_1+1367M_1$&$559$&$I_1+23501M_1$\\
			(This work)&&$+12317S_1$&&\\
			\hline
		\end{tabular}
		\caption{Comparison of the cost of the  Miller loop and the final exponentiation at $128$-bit security level.}\label{compare128}
	\end{table}
\end{center}

%%%%%%%%%%%%%%%%%%%%%%%%%%%
\subsection{New parameters and costs for pairings at the 192-bit security levels for $k=15$ and $k=27$}

\begin{itemize}
	\item Case of $k=15$. Following the recommendation from Table \ref{RNP}, we found the value $x=2^{72}+2^{40}+2^{9}+2^{5}+1$. This gives a prime $p$ of 863 bits and a prime $r$ of 577 bits. We proceed as described in Section \ref{k=15} to obtain the cost of the Miller loop and the final exponentiation. The Miller loop in this case costs $4(15M_1+13M_5+3S_5)+72(15M_1+6M_5+7S_5)+71S_{15}+75M_{15}$. This is equal to $75M_{15}+484M_5+1140M_1+71S_{15}+516S_5$. Using the arithmetic in Table \ref{COSTARITHMETICff}, the overall cost is $8871M_1+7839S_1$.  Using the value of $x$ given above,  the hard part of the final exponentiation costs  $2(72S_{15}+3M_{15})+9(72S_{15}+4M_{15})+20M_{15}+1S_{15}+4I_{G_{\phi_3(p^5)}}=62M_{15}+793S_{15}+4I_{G_{\phi_3(p^5)}}$ and $p,p^2, p^3,p^4,p^5,p^6, p^7, p^8, p^9$-Frobenius maps. The total cost of the final exponentiation is $ I_{15}+63M_{15}+793S_{15}+4I_{G_{\phi_3(p^5)}}$ and $p,p^2, p^3,p^4,2*p^5,p^6, p^7, p^8, p^9$-Frobenius maps for a total cost of $ I_1+3201M_1+35816S_1$.
	\item Case of $k=27$. Following the recommendation from Table \ref{RNP}, we found the value $x=2^{25}+2^{14}+2^{17}+2^{4}+1$. This gives a prime $p$ of 511 bits and a prime factor of $r$ of 410 bits. We proceed as described in Section \ref{k=27} to obtain the cost of the Miller loop and the final exponentiation. The Miller loop in this case costs $4(9M_1+1I_9+2S_9+3M_9)+25(9M_1+1I_9+2S_9+3M_9)+24S_{27}+27M_{27}$. This is equal to $ 29I_{9}+27M_{27}+87M_9+261M_1+24S_{27}+58S_9$. Using the arithmetic in Table \ref{COSTARITHMETICff}, the overall cost is $29I_1+11052M_1+4798S_1$.  Using the value of $x$ given above,  the final exponentiation costs  $I_{27}+2(25S_{27}+3M_{27})+17(25S_{27}+4M_{27})+11M_{27}+2I_{G_{\phi_3(p^9)}}$ and $p,p^2, p^3,p^4,p^5,p^6, p^7, p^8, 2*p^9$-Frobenius maps for a total cost of the final exponentiation  $ I_1+19191M_1+59587S_1$.
	
	The cost for the case $k=24$ are obtained with the parameter given in \cite{DuqBar18} and the formulas from \cite{AraFueMen12}
\end{itemize}

Table \ref{compare192} below compares our results with previous results at the $192$-bit security level.
\begin{center}
	\footnotesize
	\begin{table}[h]
		\begin{tabular}{|c|c|c|c|c|}
			
			\hline
			\textbf{Curves-Ref.}  & \textbf{Miller loop} &\textbf{Final Exp.} &Size of $p$&Total($S_1=M_1$)\\
			
			\hline
			$k=15$ & $8871M_1+7839S_1$ &$ I_1+3201M_1$&$863$&$I_1+55727M_1$\\
			(This work) &&$+35816S_1$&&\\
			\hline
			BLS$27$&$29I_1+11052M_1$&$ I_1+19191M_1$&$511$&$30I_1+94628M_1$\\
			
			(This work)&$+4798S_1$&$+59587S_1$&&\\
			
			\hline
			KSS$18$ \cite{DuqBar18}&$15270M_1+2590S_1$&$8I_1+7977M_1$&$677$&$8I_1+44147M_1$\\
			&&$+18310S_1$&&\\
			\hline
			BLS$24$ \cite{DuqBar18,AraFueMen12}&$15495M_1$&$10I_1+27914M_1$&$554$&$10I_1+43409M_1$\\
			\hline
		\end{tabular}
		\caption{Comparison of the cost of the  Miller loop and the final exponentiation at $192$-bit security level.}\label{compare192}
	\end{table}
\end{center}

%%%%%%%%%%%%%%%%%%%%%%%
\subsection{New parameters and costs for pairings at the 256-bit security levels for $k=27$ and $k=24$}

\begin{itemize}
	\item Case of $k=27$. Following the recommendation from Table \ref{RNP}, we found the value $x=2^{51}+2^{42}+2^{28}+2^{9}+1$. This gives a prime $p$ of 1019 bits and a prime factor of $r$ of 883 bits. We proceed as described in Section \ref{k=27} to obtain the cost of the Miller loop and the final exponentiation. The Miller loop in this case costs $4(9M_1+1I_9+2S_9+3M_9)+51(9M_1+1I_9+2S_9+3M_9)+50S_{27}+53M_{27}$. This is equal to $55I_{9}+53M_{27}+165M_9+495M_1+50S_{27}+110S_9$. Using the arithmetic in Table \ref{COSTARITHMETICff}, the overall cost is $55I_1+21348M_1+9660S_1$.  Using the value of $x$ given above,  the final exponentiation costs  $I_{27}+2(51S_{27}+3M_{27})+17(51S_{27}+4M_{27})+11M_{27}+2I_{G_{\phi_3(p^9)}}$ and $p,p^2, p^3,p^4,p^5,p^6, p^7, p^8, 2*p^9$-Frobenius maps for a total cost of the final exponentiation  $ I_1+19191M_1+122337S_1$.
	
	The cost for the case $k=24$ is obtained with the parameter given in \cite{DuqBar18} and the formulas from \cite{AraFueMen12}
\end{itemize}

Table \ref{compare256} compares our results with previous results at the 256-bit security level.
\begin{center}
	\footnotesize
	\begin{table}[h]
		\begin{tabular}{|c|c|c|c|c|}
			
			\hline
			\textbf{Curves-Ref.}  & \textbf{Miller loop} &\textbf{Final Exp.} &Size of $p$&Total($S_1=M_1$)\\
			
			\hline
			$k=27$  & $55I_1+21348M_1$ &$ I_1+19191M_1$&$1019$&$56I_1+172536M_1$\\
			(This work) &$+9660S_1$&$+122337S_1$&&\\
			\hline
			BLS$24$\cite{DuqBar18,AraFueMen12} &$18812M_1$&$10I_1+43142M_1$&$1029$&$10I_1+61954M_1$\\
			\hline
			KSS$18$\cite{DuqBar18}&$32238M_1+2620S_1$&$8I_1+7977M_1$&$1495$&$8I_1+82355M_1$\\
			&&$+39520S_1$&&\\
			\hline
			BLS$48$\cite{FOUOBAN}&$34778M_1$&$19I_1+110212M_1$&$575$&$19I_1+144990M_1$\\
			\hline
			Aur.$54$\cite{FOUOBAN}&$39976M_1+5200S_1$&$27I_1+256147M_1$&$569$&$27I_1+313397M_1$\\
			&&$+12074S_1$&&\\
			\hline
		\end{tabular}
		\caption{Comparison of the cost of the  Miller loop and the final exponentiation at the $256$-bit security level.}\label{compare256}
	\end{table}
\end{center}

%%%%%%%%%%%%%%%%%%%%%%%
\subsection{Comparison}\label{compare}

To make a fair comparison of the results in Tables \ref{compare128}, \ref{compare192} and \ref{compare256}, we take note of the size of the base field. We consider implementations on a 64-bit platform. Then, following \cite{AraFueMen12}, an $\F_{p}$-element is represented with $\ell=1+  \mbox{log}_2(p) $ binary coefficients packed in $n_{64}=\lceil \frac{\ell}{64}\rceil 64$ bits processor words and an $\F_{p}$-multiplication can be implemented with approximately $2n_{64}^{2}+n_{64}$ operations. We let $m_c$ denote the cost of a multiplication in the finite field $\F_{p}$ where $p$ is of $c$ bits. For Table \ref{compare128} we have that $m_{461}\approx 1.35m_{371}$. From this we see that, at the 128-bit security level, the total cost of computing the optimal ate pairing for elliptic curves with $k=15$ is $19199m_{461}$ making these curves faster than the well known BN curves but slower than the KSS16 curves found in \cite{DuqBar18} as the best one at the 128-bit security level.  From Table \ref{compare192} we have that the cost for computing the optimal ate pairing for curves with $k=15$ is $166067m_{511}$  as $m_{863}\approx 2.98m_{511}$ . We conclude that at the 192-bit security level  computing the optimal ate pairing is faster on elliptic curves with embedding degree $k=27$ than on curves with $k=15$ and in this case the BLS24 curves remain faster. At the 256-bit security level, we have that the BLS24 curves are the faster.

%%%%%%%%%%%%%%%%%%%%%%%%%
\section{Conclusion}\label{conclusion}

In this work we provided details and important improvements in the computation of the Miller loop and the final exponentiation for the optimal ate pairing on elliptic curves admitting cubic twists. An explicit cost evaluation is given for the Miller loop in the case of elliptic curves of embedding degrees 9 and 15. In particular for $k=15$ and $k=27$, we obtained an improvement, in terms of operations in the base field, of up to 25\% and 29\% respectively in the computation of the final exponentiation. We also obtained that elliptic curves with embedding degree $k=15$ present faster results than BN12  curves at the 128-bit security levels.  In comparison with curves having even embedding degrees we find that a lot of improvements are still required in the case of curves with odd embedding degree. One could try to compute compressed squarings in this case.  To ensure the correctness of all the formulas used in this work, a MAGMA code for the implementation of the three pairings is given. Furthermore, a brief look at the parameters used in this work reveals that the curves with odd embedding degrees studied in this work are not \textit{subgroup secure} ordinary curves \cite{BarCosMisNae15} and are not protected against \textit{small-subgroup attacks} \cite{LimLeeL97}. However this is not a particular case of elliptic curves of odd embedding degree but it appears from \cite{BarCosMisNae15} that most of  such parameters that have been found for curves with even embedding degree such as BN12 curves \cite{BarNar05}, KSS16
curves \cite{kachisa08constructing} or BLS12 curves \cite{BarLynSco02}; do not satisfied these security properties. In future work we could search for parameters to fulfill this security issue.

%%%%%%%%%%%%%%%%%%%%%%%%%%
\section*{Acknowledgment}

The authors thank the anonymous reviewers who helped in improving the quality of this work. A particular thank goes to the editorial board and the linguistic editor of the Journal of Groups, Complexity and Cryptology for the relevant linguistic corrections which significantly improved this work.

%%%%%%%%%%%%%%%%%%%%%%%%%%
%%%%%%%%%%%%%%%%%%%%%%%%%%
\bibliographystyle{plain}
\bibliography{bibliographyEAN}

%%%%%%%%%%%%%%%%%%%%%%%%%%
%%%%%%%%%%%%%%%%%%%%%%%%%%
\appendix

%%%%%%%%%%%%%%%%%%%%
\section{Arithmetic in $\F_{p^{9}}$}\label{arithFP9}
%\subsection{Squaring in $\F_{p^{9}}$ }

Let $a= a_{0}+ a_{1}v+ a_{2}v^{2} \in \mathbb{F}_{p^{9}}$ with $ a_{i}\in\F_{p^{3}}$.
%$a^2= A_{0}+A_{1}v+A_{2}v^{2}$ where $$ \left\{ \begin{array}{ccc}
%A_{0} & = & a_{0}^{2}+2a_{1}a_{2}7^{1/3} \\
%A_{1} & = & 2a_{0}a_{1}+a_{2}^{2}7^{1/3} \\
%A_{2} & = & 2a_{0}a_{2} +a_{1}^{2}
%\end{array} \right. $$ \\
%This costs $3M_3+3S_3$ or $6S_3$ since the computation of $2xy$ can be done as $(x+y)^2-x^2-y^2$ when the squares $x^2$ and $y^2$ are known.

%%%%%%%%%%%%%%%%%%%%
\subsection{Cyclotomic inversion}\label{ci9}

We assume that $a$ lies in the cyclotomic subgroup  $G_{\phi_{3}(p^{3})}$, so that $a^{p^{6}+p^{3}+1}=1$, i.e., $a^{-1}=a^{p^{6}}a^{p^{3}}$.
In order to compute $a^{p^{6}}a^{p^{3}}$, we need the values of $v^{p^{3}}$ and $v^{p^{6}}$. But $v^{p^{3}} = v^{3(p^{3}-1)/3+1} = v^{3(p^{3}-1)/3}v = (v^{3})^{(p^{3}-1)/3}v = (7^{1/3})^{(p^{3}-1)/3}v$ since $v^{3}=7^{1/3}$.

Let $\mu =(7^{1/3})^{(p^{3}-1)/3}$; we have $\mu \neq 1$  and $\mu^{3}=1$ so that  $\mu$ is a primitive cubic root of unity in $\mathbb{F}_{p^{3}}$.
We obtain  $v^{p^{3}}=\mu v$ and $v^{p^{6}}=(v^{p^{3}})^{p^{3}}= (\mu v)^{p^{3}}= \mu (v)^{p^{3}}=\mu \mu v=  \mu ^{2}v$.
We then have $a^{p^{3}}= a_{0}^{p^{3}}+ a_{1}^{p^{3}}v^{p^{3}}+ a_{2}^{p^{3}}(v^{2})^{p^{3}}=a_{0}+ a_{1}v^{p^{3}}+ a_{2}(v^{2})^{p^{3}}=a_{0}+ a_{1}\mu v+ a_{2}\mu^{2}v^{2}$ and $a^{p^{6}}= (a^{p^{3}})^{p^{3}}=a_{0}+ a_{1}(\mu v)^{p^{3}}+ a_{2}(\mu^{2}v^{2})^{p^{3}}=a_{0}+ a_{1}\mu^{2} v+ a_{2}\mu^{4}v^{2} $. So, when using $v^{3}=7^{1/3}$ and $\phi_{3}(\mu)= \mu^{2}+\mu +1=0$,
we finally have:
$$a^{p^{6}}a^{p^{3}}= (a_{0}^{2}-a_{1}a_{2}7^{1/3})+(a_{2}^{2}7^{1/3}-a_{0}a_{1})v+
(a_{1}^{2}-a_{0}a_{2}) v^{2}$$
This costs $3M_3+3S_3=18M_1+15S_1$ with additional additions.

%%%%%%%%%%%%%%%%%%%%%%%%%%%%
\subsection{Frobenius operators}

The $p^{i}$-Frobenius is the map  $\pi^{i}\colon \F_{p^{9}}  \longrightarrow   \F_{p^{9}}$ given by $\pi^{i}(a) = a^{p^{i}}$.
Let $a \in \F_{p^{9}}$,  $a= a_{0}+ a_{1}v+ a_{2}v^{2}$ with $ a_{i}\in\mathbb{F}_{p^{3}}$; then $\pi(a)=a_{0}^{p}+a_{1}^{p}v^{p}+a_{2}^{p}(v^{2})^{p}$. Now $a_{0}\in \mathbb{F}_{p^{3}}$ can be written as $a_{0}= g_{0}+g_{1}u+g_{2}u^{2}$, $g_{i}\in \mathbb{F}_{p}$ so that $a_{0}^{p}= g_{0}+g_{1}u^{p}+g_{2}(u^{2})^{p}$.

We have $u^{p}=u^{3(p-1)/3+1}=(u^{3})^{(p-1)/3}u=7^{(p-1)/3}u$ and; since $7$ is not a cube in $ \mathbb{F}_{p}$, $7^{(p-1)/3} \neq 1$.
Let $\alpha = 7^{(p-1)/3}$; then $\alpha\neq 1$  and $\alpha^{3}=1$. This means that $\alpha$ is a primitive cubic root of unity in $\mathbb{F}_{p}$ and  $u^{p}=\alpha u$. Therefore $a_{0}^{p}= g_{0}+g_{1}u^{p}+g_{2}(u^{2})^{p}=g_{0}+g_{1}\alpha u+g_{2}\alpha^{2} u^{2}$ and; similarly;
$a_{1}^{p}= g_{3}+g_{4}u^{p}+g_{5}(u^{2})^{p}=g_{3}+g_{4}\alpha u+g_{5}\alpha^{2} u^{2}$ and
$a_{2}^{p}= g_{6}+g_{7}u^{p}+g_{8}(u^{2})^{p}=g_{6}+g_{7}\alpha u+g_{8}\alpha^{2} u^{2}$. Now, for the computation of $v^p$, observe that
$v^{p}=v^{3(p-1)/3+1}=(v^{3})^{(p-1)/3}v=(7^{1/3})^{(p-1)/3}v=7^{(p-1)/9}v$ so that, if $\beta=7^{(p-1)/9}$, then
we have $\beta \neq 1$, $\beta^{3}=7^{(p-1)/3}=\alpha\neq 1$, $\beta^{9}=1$.
Thus $\beta$ is a primitive ninth root of unity in $\mathbb{F}_{p}$ and $v^{p}=\beta v$.

Finally $a^{p}=g_{0}+g_{1}\alpha u+g_{2}\alpha^{2}u^{2}+(g_{3}\beta+g_{4}\alpha \beta u+g_{5}\alpha^{2}\beta u^{2})v+(g_{6}\beta^{2}+g_{7}\alpha \beta^{2}u+g_{8}\alpha^{2}\beta^{2}u^{2})v^{2}$ and
the following algebraic relations: $\alpha = \beta^{3}$, $\alpha\beta=\beta^{4}$, $\alpha \beta^{2}=\beta^{5}$, $\alpha^{2}\beta=\beta^{7}$, $\alpha^{2}\beta^{2}=\beta^{8}$ yield to $a^{p}=(g_{0}+g_{1}\beta^{3}u+g_{2}\beta^{6}u^{2})+(g_{3}\beta+g_{4}\beta^{4}u+g_{5}\beta^{7}u^{2})v+(g_{6}\beta^{2}+g_{7}\beta^{5}u
+g_{8}\beta^{8}u^{2})v^{2}$.
The cost of $p$-Frobenius is $8M_{1}$. This is the same as the cost of $p^2,p^4,p^5,p^7$ and $p^8$-Frobenius.
For the $p^3$-Frobenius operator, observe from \ref{ci9} that  $v^{p^{3}}=\mu v$. Then
$$a^{p^{3}}= a_{0}+ a_{1}\mu v+ a_{2}\mu^{2}v^{2}=(g_{0}+g_{1}u+g_{2}u^{2})
+(g_{3}+g_{4}u+g_{5}u^{2})\mu v +(g_{6}+g_{7}u+g_{8}u^{2}) \mu^{2}v^{2}.$$
As $t= \mu^{2}$ is precomputed; we finally have
$$a^{p^{3}}=(g_{0}+g_{1}u+g_{2}u^{2}) +(g_{3}\mu+g_{4}\mu u+g_{5}\mu u^{2}) v +(g_{6}t+g_{7} t u +g_{8}tu^{2}) v^{2}.$$
The cost of $p^{3}$-Frobenius: $6M_1$. This is the same as the cost of $p^6$-Frobenius.

%%%%%%%%%%%%%%%%%%%%%%%%%%%%
%%%%%%%%%%%%%%%%%%%%%%%%%%%%
\section{Arithmetic in  $\F_{p^{27}}$ }\label{arithFP27}

%%%%%%%%%%%%%%%%%%%%%%%%%%%%
\subsection{Cyclotomic inversion}\label{ci27}

We follow the same procedure as in \ref{ci9}.
The element $a= a_{0}+ a_{1}w+ a_{2}w^{2} \in \mathbb{F}_{p^{27}}$ with $ a_{i}\in\mathbb{F}_{p^{9}}$ in the cyclotomic subgroup  $G_{\phi_{3}(\mathbb{F}_{p^{9}})}$ satisfies  $a^{p^{18}+p^{9}+1}=1$ so that  $a^{-1}=a^{p^{18}}a^{p^{9}}$.

In order to compute $a^{p^{18}}a^{p^{9}}$, we need the values of $w^{p^{9}}$ and $w^{p^{18}}$. We have
$$w^{p^{9}} = w^{3(p^{9}-1)/3+1} = w^{3(p^{9}-1)/3}w
= (w^{3})^{(p^{9}-1)/3}w = (7^{1/9})^{(p^{9}-1)/3}w$$
since $w^{3}=7^{1/9}$.
Let $\sigma =(7^{1/9})^{(p^{9}-1)/3}$; then $\sigma \neq 1$  and $\sigma^{3}=1$. Hence $\sigma$ is a primitive cubic root of unity in $\mathbb{F}_{p^{9}}$, i.e. $\phi_{3}(\sigma)=0$. We obtain  $w^{p^{9}}=\sigma w$ and we  now compute  $w^{p^{18}}$ as
\begin{align*}
w^{p^{18}}&=(w^{p^{9}})^{p^{9}}= (\sigma w)^{p^{9}}= \sigma (w)^{p^{9}}=\sigma \sigma w=  \sigma ^{2}w,\\
a^{p^{9}}& =a_{0}+ a_{1}w^{p^{9}}+ a_{2}(w^{2})^{p^{9}}=a_{0}+ a_{1}\sigma w+ a_{2}\sigma^{2}w^{2},\text{ and} \\
a^{p^{18}}&= (a^{p^{9}})^{p^{9}}=a_{0}+ a_{1}(\sigma w)^{p^{9}}+ a_{2}(\sigma^{2}w^{2})^{p^{9}}=a_{0}+ a_{1}\sigma^{2} w+ a_{2}\sigma^{4}w^{2}.
\end{align*}
%Finally
%$a^{p^{18}}a^{p^{9}} = a_{0}^{2}+a_{0}a_{1}\sigma^{2}w+  a_{0}a_{2}\sigma w^{2}+ a_{0}a_{1}\sigma w+a_{1}^{2}w^{2}+a_{1}a_{2}\sigma^{2}7^{1/9}+a_{0}a_{2}\sigma^{2}w^{2}+
% a_{1}a_{2}\sigma 7^{1/9}+a_{2}^{2}7^{1/9}w$.\\
%
After expanding and reducing, using $w^{3}=7^{1/9}$ and $\phi_{3}(\sigma)= \sigma^{2}+\sigma +1=0$, we obtain
$$a^{p^{18}}a^{p^{9}}= (a_{0}^{2}-a_{1}a_{2}7^{1/9})+(a_{2}^{2}7^{1/9}-a_{0}a_{1})w+
(a_{1}^{2}-a_{0}a_{2}) w^{2}.$$
The computation costs $3(36M_1)+3(25S_1)=108M_1+75S_1$.

%%%%%%%%%%%%%%%%%%%%%%%%%%%%
\subsection{Frobenius operators}

The $p^{i}-$Frobenius is the map  $\pi^{i}\colon: \F_{p^{27}}  \longrightarrow   \F_{p^{27}}$ given by $\pi^i(a) = a^{p^{i}}$.
Let $a= a_{0}+ a_{1}w+ a_{2}w^{2}$ with $ a_{i}\in\mathbb{F}_{p^{9}}$ an element of $ \mathbb{F}_{p^{27}}$. Then
$$\pi(a)=a^{p}=(a_{0}+ a_{1}w+ a_{2}w^{2})^{p}=a_{0}^{p}+a_{1}^{p}w^{p}+a_{2}^{p}(w^{2})^{p}.$$
The element $a_{0}\in \mathbb{F}_{p^{9}}$ can be written as
\begin{align*}
a_{0} &= (h_{0}+h_{1}u+h_{2}u^{2})+(h_{3}+h_{4}u+h_{5}u^{2})v+ (h_{6}+h_{7}u+h_{8}u^{2})v^{2}, \text{ $h_{i}\in \F_{p}$ and then}, \\
a_{0}^{p} &= (h_{0}+h_{1}u+h_{2}u^{2}+(h_{3}+h_{4}u+h_{5}u^{2})v+ (h_{6}+h_{7}u+h_{8}u^{2})v^{2})^{p}, \\
h_{i}^{p} &=h_{i}, \\
u^{p} &=u^{3(p-1)/3+1}=(u^{3})^{(p-1)/3}u=7^{(p-1)/3}u.
\end{align*}
Since $7$ is not a cube in $\F_{p}$, we have $\alpha = 7^{(p-1)/3}$ $\alpha\neq 1$ and $\alpha^{3}=1$. It means that $\alpha$ is a primitive cubic root of unity in $\mathbb{F}_{p}$ and $u^{p}=\alpha u$. $v^{p}=v^{3(p-1)/3+1}=(v^{3})^{(p-1)/3}v=(7^{1/3})^{(p-1)/3}v=7^{(p-1)/9}v$.

We have $\beta=7^{(p-1)/9}\neq 1$ and $\beta^{9}=1$. Thus $\beta$ is a primitive ninth root of unity in $\mathbb{F}_{p}$ and $v^{p}=\beta v$. Also $w^{p}=w^{3(p-1)/3+1}=(w^{3})^{(p-1)/3}v=(7^{1/9})^{(p-1)/3}v=7^{(p-1)/27}v$.
We also observe that $\gamma=7^{(p-1)/27} \neq 1$, $\gamma^{3}=7^{(p-1)/9}=\beta \neq 1$, $\gamma^{9}=7^{(p-1)/3}=\alpha \neq 1$, $\gamma^{27}=1$.
Thus $\gamma$ is a primitive twenty-seventh root of unity in $\F_{p}$ and $w^{p}= \gamma w$.
\begin{align*}
a_{0}^{p} &=((h_{0}+h_{1}u+h_{2}u^{2})+(h_{3}+h_{4}u+h_{5}u^{2})v+(h_{6}+h_{7}u+h_{8}u^{2})v^{2})^{p} \\
&=(h_{0}+h_{1}u^{p}+h_{2}(u^{2})^{p})+(h_{3}+h_{4}u^{p}+h_{5}(u^{2})^{p})v^{p}+ (h_{6}+h_{7}u^{p}+h_{8}(u^{2})^{p})(v^{2})^{p} \\
&=(h_{0}+h_{1}\alpha u+h_{2}\alpha^{2} u^{2})+(h_{3}+h_{4} \alpha u +h_{5} \alpha^{2} u^{2})\beta v+ (h_{6}+h_{7}\alpha u+h_{8}\alpha^{2} u^{2}) \beta^{2}v^{2} \\
&=(h_{0}+h_{1}\alpha u+h_{2}\alpha^{2} u^{2})+(h_{3} \beta +h_{4} \alpha \beta u +h_{5} \alpha^{2}\beta u^{2}) v \\
&\hskip .75cm + (h_{6} \beta^{2} +h_{7}\alpha \beta^{2} u+h_{8}\alpha^{2} \beta^{2}u^{2}) v^{2},
\end{align*}
\begin{align*}
a_{1}^{p} &= (h_{9}+h_{10}u+h_{11}u^{2})+(h_{12}+h_{13}u+h_{14}u^{2})v+ (h_{15}+h_{16}u+h_{17}u^{2})v^{2})^{p} \\
&=(h_{9}+h_{10}u^{p}+h_{11}(u^{2})^{p})+(h_{12}+ h_{13}u^{p}+h_{14}(u^{2})^{p})v^{p}\\
&\hskip .75cm + (h_{15}+h_{16}u^{p}+h_{17}(u^{2})^{p})(v^{2})^{p} \\
&=(h_{9}+h_{10}\alpha u+h_{11}\alpha^{2} u^{2})+(h_{12}+h_{13} \alpha u +h_{14} \alpha^{2} u^{2})\beta v\\
&\hskip .75cm + (h_{15}+h_{16}\alpha u+h_{17}\alpha^{2} u^{2}) \beta^{2}v^{2} \\
&=(h_{9}+h_{10}\alpha u+h_{11}\alpha^{2} u^{2})+(h_{12} \beta +h_{13} \alpha \beta u +h_{14} \alpha^{2}\beta u^{2}) v \\
&\hskip .75cm  + (h_{15} \beta^{2} +h_{16}\alpha  \beta^{2}u+h_{17}\alpha^{2} \beta^{2}u^{2}) v^{2},
\end{align*}
\begin{align*}
a_{2}^{p} &= (h_{18}+h_{19}u+h_{20}u^{2})+(h_{21}+h_{22}u+h_{23}u^{2})v+ (h_{24}+h_{25}u+h_{26}u^{2})v^{2})^{p} \\
&=(h_{18}+h_{19}u^{p}+h_{20}(u^{2})^{p})+(h_{21}+
h_{22}u^{p}+h_{23}(u^{2})^{p})v^{p}\\
&\hskip .75cm + (h_{24}+h_{25}u^{p}+h_{26}(u^{2})^{p})(v^{2})^{p} \\
&=(h_{18}+h_{19}\alpha u+h_{20}\alpha^{2} u^{2})+(h_{21}+h_{22} \alpha u +h_{23} \alpha^{2} u^{2})\beta v\\
&\hskip .75cm + (h_{24}+h_{25}\alpha u+h_{26}\alpha^{2} u^{2}) \beta^{2}v^{2} \\
&=(h_{18}+h_{19}\alpha u+h_{20}\alpha^{2} u^{2})+(h_{21} \beta +h_{22} \alpha \beta u +h_{23} \alpha^{2}\beta u^{2}) v \\
&\hskip .75cm + (h_{24} \beta^{2} +h_{25}\alpha \beta^{2} u+ h_{26}\alpha^{2} \beta^{2}u^{2}) v^{2},
\end{align*}
\begin{align*}
\pi(a) &=(a_{0}+ a_{1}w+ a_{2}w^{2})^{p} =a_{0}^{p}+a_{1}^{p}w^{p}+a_{2}^{p}(w^{2})^{p} =a_{0}^{p}+a_{1}^{p}\gamma w+a_{2}^{p}\gamma^{2}w^{2} \\
&=(h_{0}+h_{1}\alpha u+h_{2}\alpha^{2} u^{2})+(h_{3} \beta +h_{4} \alpha \beta u +h_{5} \alpha^{2}\beta u^{2}) v \\
&\hskip .75cm + (h_{6} \beta^{2} +h_{7}\alpha \beta^{2} u+h_{8}\alpha^{2} \beta^{2}u^{2}) v^{2} \\
&\hskip .75 cm +((h_{9}+h_{10}\alpha u+h_{11}\alpha^{2} u^{2})+(h_{12} \beta +h_{13} \alpha \beta u +h_{14} \alpha^{2}\beta u^{2}) v \\
&\hskip .75cm + (h_{15} \beta^{2} +h_{16}\alpha \beta^{2} u+ h_{17}\alpha^{2} \beta^{2}u^{2}) v^{2})\gamma w+ ((h_{18}+h_{19}\alpha u+h_{20}\alpha^{2} u^{2})\\
&\hskip .75cm +(h_{21} \beta +h_{22} \alpha \beta u +h_{23} \alpha^{2}\beta u^{2}) v +(h_{24} \beta^{2} +h_{25}\alpha \beta^{2} u+h_{26}\alpha^{2} \beta^{2}u^{2}) v^{2}) \gamma^{2}w^{2}.
\end{align*}
We have the following algebraic relations:
$\alpha = \beta^{3}$, $\alpha\beta=\beta^{4}$, $\alpha \beta^{2}=\beta^{5}$, $\alpha^{2}\beta=\beta^{7}$ and $\alpha^{2}\beta^{2}=\beta^{8}$.
Therefore
\begin{align*}
\pi(a) &=((h_{0}+h_{1}\beta^{3} u+h_{2}\beta^{6} u^{2})+(h_{3} \beta +h_{4}  \beta^{4} u +h_{5} \beta^{7} u^{2})v + (h_{6} \beta^{2} +h_{7}\beta^{5} u +h_{8} \beta^{8}u^{2})v^{2}) \\
&\hskip .75cm +((h_{9}\gamma+h_{10}\beta^{3} \gamma u+h_{11}\beta^{6} \gamma u^{2})+(h_{12} \beta \gamma +h_{13} \beta^{4} \gamma u +h_{14} \beta^{7} \gamma u^{2})v \\
&\hskip .75cm  + (h_{15} \beta^{2}\gamma +h_{16}\beta^{5}\gamma u+h_{17}\beta^{8}\gamma u^{2}) v^{2})w+ ((h_{18}\gamma^{2}+h_{19}\beta^{3}\gamma^{2} u+h_{20}\beta^{6}\gamma^{2} u^{2})  \\
&\hskip .75cm +(h_{21} \beta \gamma^{2}  +h_{22} \beta^{4} \gamma^{2}u +h_{23} \beta^{7}\gamma^{2} u^{2})v + (h_{24} \beta^{2}\gamma^{2} +h_{25} \beta^{5}\gamma^{2} u + h_{26} \beta^{8}\gamma^{2}u^{2})v^{2})w^{2}.
\end{align*}
The following values are precomputed:
$\lambda_{0}=\beta^{2}$, $\lambda_{1}=\beta^{3}$, $\lambda_{2}=\beta^{4}$,
$\lambda_{3}=\beta^{5}$, $\lambda_{4}=\beta^{6}$,
$\lambda_{5}=\beta^{7}$, $\lambda_{6}=\beta_{8}$, $\lambda_{7}= \gamma^{2}$
, $\lambda_{8}=\beta \gamma$,
$\lambda_{9}=\lambda_{0}\gamma$,
$\lambda_{10}=\lambda_{1}\gamma$, $\lambda_{11}=\lambda_{2}\gamma$, $\lambda_{12}=\lambda_{3}\gamma$,
$\lambda_{13}=\lambda_{4}\gamma $, $\lambda_{14}=\lambda_{5}\gamma $,
$\lambda_{15}=\lambda_{6}\gamma$,
$\lambda_{16}=\lambda_{0}\lambda_{7}$, $\lambda_{17}=\lambda_{1}\lambda_{7}$,
$\lambda_{18}=\lambda_{2}\lambda_{7}$, $\lambda_{19}=\lambda_{3}\lambda_{7}$,
$\lambda_{20}=\lambda_{4}\lambda_{7}$, $\lambda_{21}=\lambda_{5}\lambda_{7}$,
$\lambda_{22}=\lambda_{6}\lambda_{7}$. $\lambda_{23}=\beta \lambda_{7}$.
Thus
\begin{align*}
\pi(a) &= ((h_{0}+h_{1}\lambda_{1} u+h_{2}\lambda_{4} u^{2})+(h_{3} \beta +h_{4}\lambda_{2} u +h_{5} \lambda_{5} u^{2})v \\
&\hskip .75cm  + (h_{6} \lambda_{0} +h_{7}\lambda_{3}u + h_{8} \lambda_{6}u^{2})v^{2}) + ((h_{9}\gamma+h_{10}\lambda_{10} u+h_{11}\lambda_{13} u^{2}) \\
&\hskip .75cm +(h_{12} \lambda_{8} +h_{13} \lambda_{11} u +h_{14} \lambda_{14} u^{2})v + (h_{15} \lambda_{9} + h_{16}\lambda_{12} u+h_{17}\lambda_{15} u^{2}) v^{2})w \\
&\hskip .75cm +((h_{18}\lambda_{7}+h_{19}\lambda_{17}u+h_{20}\lambda_{20} u^{2})+(h_{21} \lambda_{23}  + h_{22} \lambda_{18}u + h_{23} \lambda_{21} u^{2})v \\
&\hskip .75cm  + (h_{24} \lambda_{16} +h_{25} \lambda_{19} u + h_{26} \lambda_{22}u^{2})v^{2})w^{2}.
\end{align*}
The cost of $p$-Frobenius is $26M_1+18a$. This is also equal to the cost of $p^2,p^4, p^5, p^7, p^8$-Frobenius.
For the $p^9$-Frobenius operator, observe from \ref{ci27} that  $w^{p^{9}}=\sigma w$. Then
%$a^{p^{9}}=(a_{0}+ a_{1}w+ a_{2}w^{2})^{p^{9}}= a_{0}^{p^{9}}+ a_{1}^{p^{9}}w^{p^{9}}+ a_{2}^{p^{9}}(w^{2})^{p^{9}}=a_{0}+ a_{1}w^{p^{9}}+ a_{2}(w^{2})^{p^{9}}$
%
%$\;\ \;\ \;\ \;\ =a_{0}+ a_{1}\sigma w+ a_{2}\sigma^{2}w^{2}$. \\
%
%$a_{0}= (h_{0}+h_{1}u+h_{2}u^{2})+(h_{3}+h_{4}u+h_{5}u^{2})v+ (h_{6}+h_{7}u+h_{8}u^{2})v^{2}$, \\
%
%$a_{1}= (h_{9}+h_{10}u+h_{11}u^{2})+(h_{12}+h_{13}u+h_{14}u^{2})v+ (h_{15}+h_{16}u+h_{17}u^{2})v^{2}$, \\
%
%$a_{2}= (h_{18}+h_{19}u+h_{20}u^{2})+(h_{21}+h_{22}u+h_{23}u^{2})v+ (h_{24}+h_{25}u+h_{26}u^{2})v^{2}$, \;\ $h_{i}\in \mathbb{F}_{p}$.\\
%
\begin{align*}
a^{p^{9}} &= a_{0}+ a_{1}\sigma w+ a_{2}\sigma^{2}w^{2} \\
&=((h_{0}+h_{1}u+h_{2}u^{2})+(h_{3}+h_{4}u+h_{5}u^{2})v+ (h_{6}+h_{7}u + h_{8}u^{2})v^{2}) \\
&\hskip .75cm +((h_{9}+h_{10}u+h_{11}u^{2})+(h_{12}+h_{13}u+h_{14}u^{2})v+ (h_{15}+h_{16}u+h_{17}u^{2})v^{2})\sigma w \\
&\hskip .75cm + ((h_{18}+h_{19}u+h_{20}u^{2})+(h_{21}+h_{22}u+h_{23}u^{2})v+ (h_{24}+h_{25}u+h_{26}u^{2})v^{2})\sigma^{2}w^{2}.
\end{align*}
We then have
\begin{align*}
a^{p^{9}} &= ((h_{0}+h_{1}u+h_{2}u^{2})+(h_{3}+h_{4}u+h_{5}u^{2})v+ (h_{6}+h_{7}u+h_{8}u^{2})v^{2}) \\
&\hskip .75cm + ((h_{9}\sigma +h_{10}\sigma u+h_{11}\sigma u^{2})+(h_{12}\sigma +h_{13}\sigma u+h_{14} \sigma u^{2})v \\
&\hskip .75cm + (h_{15}\sigma+h_{16}\sigma u+ h_{17}\sigma u^{2})v^{2})w + ((h_{18}\sigma^{2} +h_{19}\sigma^{2}  u+h_{20}\sigma^{2} u^{2}) \\
&\hskip .75cm +(h_{21}\sigma^{2} +h_{22} \sigma^{2} u+ h_{23} \sigma^{2} u^{2})v+ (h_{24}\sigma^{2} +h_{25}\sigma^{2} u+h_{26}\sigma^{2}  u^{2})v^{2})w^{2}.
\end{align*}
As $s=\sigma^{2}$ is precomputed, we have
\begin{align*}
a^{p^{9}} &= ((h_{0}+h_{1}u+h_{2}u^{2})+(h_{3}+h_{4}u+h_{5}u^{2})v+ (h_{6}+h_{7}u+h_{8}u^{2})v^{2}) \\
&\hskip .75cm + ((h_{9}\sigma +h_{10}\sigma u+h_{11}\sigma u^{2})+(h_{12}\sigma +h_{13}\sigma u+h_{14} \sigma u^{2})v \\
&\hskip .75cm + (h_{15}\sigma+h_{16}\sigma u+h_{17}\sigma u^{2})v^{2})w +((h_{18}s +h_{19}s  u+h_{20}s u^{2}) \\
&\hskip .75cm +(h_{21}s +h_{22} s u+h_{23} s u^{2})v+ (h_{24}s +h_{25}s u+h_{26}s  u^{2})v^{2})w^{2}.
\end{align*}
The cost of $p^{9}$-Frobenius is $18M_1$. This is the same as the cost of $p^{3}$ and $p^6$ Frobenius.

%%%%%%%%%%%%%%%%%%%%%%%
\section{Arithmetic in $\F_{p^{15}}$ }\label{arithFP15}

The arithmetic of the extension field $F_{p^{5}}$ is studied in \cite{MraGuiIon09}. In this section we only consider inversion in cyclotomic subgroup and Frobenius operators.

%%%%%%%%%%%%%%%%%%%%%%
\subsection{Cyclotomic inversion}\label{ci15}

An element $a= a_{0}+ a_{1}v+ a_{2}v^{2} \in \F_{p^{15}}$ with $a_{i}\in\F_{p^{5}}$ in the cyclotomic subgroup  $G_{\phi_{3}(p^{5})}$ satisfies $a^{p^{10}+p^{5}+1}=1$ so that $a^{-1}=a^{p^{10}}a^{p^{5}}$. We have
$v^{p^{5}} = v^{5(p^{5}-1)/5+1} = v^{5(p^{5}-1)/5}v = (v^{5})^{(p^{5}-1)/5}v = (7^{1/3})^{(p^{5}-1)/5}v$ since $v^{5}=7^{1/3}$.

Let $\omega =(7^{1/3})^{(p^{5}-1)/5}$. We have $\omega \neq 1$  and $\omega^{5}=1$.
Hence $\omega$ is a primitive fifth root of unity in $\mathbb{F}_{p^{5}}$.
We obtain  $v^{p^{5}}=\omega v$ and $v^{p^{10}}=(v^{p^{5}})^{p^{5}}= (\omega v)^{p^{5}}= \omega (v)^{p^{5}}=\omega \omega v=  \omega ^{2}v$. So
\begin{align*}
a^{p^{5}} &=(a_{0}+ a_{1}v+ a_{2}v^{2})^{p^{5}}= a_{0}^{p^{5}}+ a_{1}^{p^{5}}v^{p^{5}}+ a_{2}^{p^{5}}(v^{2})^{p^{5}} \\
&= a_{0}+ a_{1}v^{p^{5}}+ a_{2}(v^{2})^{p^{5}} =a_{0}+ a_{1}\omega v+ a_{2}\omega^{2}v^{2}, \\
a^{p^{10}}& = (a^{p^{5}})^{p^{5}}=a_{0}+ a_{1}(\omega v)^{p^{5}}+ a_{2}(\omega^{2}v^{2})^{p^{5}}=a_{0}+ a_{1}\omega^{2} v+ a_{2}\omega^{4}v^{2}, \\
a^{p^{10}}a^{p^{5}}  &=  (a_{0}+ a_{1}\omega^{2} v+ a_{2}\omega^{4}v^{2})(a_{0}+ a_{1}\omega v+ a_{2}\omega^{2}v^{2}).
\end{align*}
%$ \;\ \;\ =  a_{0}^{2}+a_{0}a_{1}\omega^{2}v+  a_{0}a_{2}\omega^{4}v^{2}+a_{0}a_{1}\omega v+a_{1}^{2}\omega^{3}v^{2}+a_{1}a_{2}\omega^{5}v^{3}+a_{0}a_{2}\omega^{2}v^{2}+
% a_{1}a_{2}\omega^{4}v^{3}+a_{2}^{2}\omega^{6}v^{4}$

%$ \;\ \;\ \;\ =  a_{0}^{2}+a_{0}a_{1}\omega^{2}v+  a_{0}a_{2}\omega^{4} v^{2}+ a_{0}a_{1}\omega v+a_{1}^{2}\omega^{3} v^{2}+a_{1}a_{2}u+a_{0}a_{2}\omega^{2}v^{2}+
% a_{1}a_{2}\omega^{4} u+a_{2}^{2}\omega u v $
After expanding and reducing using $v^{3}=u$ and $\phi_{5}(\omega)=0$ we obtain
$$a^{p^{10}}a^{p^{5}}= (a_{0}^{2}+(1+\omega^{4})a_{1}a_{2}u)+\omega(a_{2}^{2}u+(1+\omega)a_{0}a_{1})v+
\omega^{2}(a_{1}^{2}\omega+(1+\omega^{2})a_{0}a_{2}) v^{2}.$$
This costs $3(9M_1)+3(9S_1)=27M_1+27S_1$.

%%%%%%%%%%%%%%%%%%%%
\subsection{Frobenius operators}

The $p^{i}-$Frobenius is the map  $\pi^{i}\colon\mathbb{F}_{p^{15}}  \rightarrow   \mathbb{F}_{p^{15}}$ given by $\pi^i(a) = a^{p^{i}}$.
Let $a \in \mathbb{F}_{p^{15}}$, say,  $a= a_{0}+ a_{1}v+ a_{2}v^{2}$ with $ a_{i}\in\mathbb{F}_{p^{5}}$.
Then
$$\pi(a)=a^{p}=(a_{0}+ a_{1}v+ a_{2}v^{2})^{p}=a_{0}^{p}+a_{1}^{p}v^{p}+a_{2}^{p}(v^{2})^{p}.$$
As $a_{0}\in \mathbb{F}_{p^{5}}$, \textit{i.e.}, $a_{0}= g_{0}+g_{1}u+g_{2}u^{2}+g_{3}u^{3}+g_{4}u^{4}$, $g_{i}\in \mathbb{F}_{p}$, we have 
$$a_{0}^{p}= (g_{0}+g_{1}u+g_{2}u^{2}+g_{3}u^{3}+g_{4}u^{4})^{p}= g_{0}+g_{1}u^{p}+g_{2}(u^{2})^{p}+g_{3}(u^{3})^{p}+g_{4}(u^{4})^{p},$$
since $g_{i}^{p}=g_{i}$.
Now we have $u^{p}=u^{5(p-1)/5+1}=(u^{5})^{(p-1)/5}u=7^{(p-1)/5}u$. As
$7$ is not a fifth power in $ \mathbb{F}_{p}$, we have  $7^{(p-1)/5} \neq 1$.
Let $\theta = 7^{(p-1)/5}$. Then $\theta \neq 1$ and $\theta^{5}=1$, that is, $\theta$ is a primitive fifth root of unity in $\mathbb{F}_{p}$ and $u^{p}=\theta u$. Then
\begin{align*}
a_{0}^{p} &= g_{0}+g_{1}u^{p}+g_{2}(u^{2})^{p}+g_{3}(u^{3})^{p}+g_{4}(u^{4})^{p}=g_{0}+g_{1}\theta u+g_{2}\theta^{2} u^{2}+g_{3}\theta^{3}u^{3}+g_{4}\theta^{4}u^{4}, \\
a_{1}^{p} &= g_{5}+g_{6}u^{p}+g_{7}(u^{2})^{p}+g_{8}(u^{3})^{p}+g_{9}(u^{4})^{p}=g_{5}+g_{6}\theta u+g_{7}\theta^{2} u^{2}+g_{8}\theta^{3}u^{3}+g_{9}\theta^{4}u^{4}, \\
a_{2}^{p} &= g_{10}+g_{11}u^{p}+g_{12}(u^{2})^{p}+g_{13}(u^{3})^{p}+g_{14}(u^{4})^{p} \\
&=g_{10}+g_{11}\theta u+g_{12}\theta^{2} u^{2}+g_{13}\theta^{3}u^{3}+g_{14}\theta^{4}u^{4}, \\
v^{p} &=v^{5(p-1)/5+1}=(v^{5})^{(p-1)/5}v=(7^{1/3})^{(p-1)/5}v=(7^{1/3})^{(p-1)/5}v.
\end{align*}
Let $\beta=(7^{1/3})^{(p-1)/5}$. Since $7^{1/3}$ is not a fifth power in $\F_{p}$, we have $\beta \neq 1$ and $\beta^{5}=1$, that is, $\beta$ is a primitive fifth root of unity in $\mathbb{F}_{p}$ and $v^{p}=\beta v$. Then
\begin{align*}
a^{p} &=(a_{0}+ a_{1}v+ a_{2}v^{2})^{p} =a_{0}^{p}+a_{1}^{p}v^{p}+a_{2}^{p}(v^{2})^{p} \\
&= (g_{0}+g_{1}\theta u+g_{2}\theta^{2} u^{2}+g_{3}\theta^{3}u^{3}+g_{4}\theta^{4}u^{4}) + (g_{5}+g_{6}\theta u+g_{7}\theta^{2} u^{2}+g_{8}\theta^{3}u^{3}+g_{9}\theta^{4}u^{4})v^{p} \\
&\hskip .75cm +(g_{10}+g_{11}\theta u+g_{12}\theta^{2} u^{2}+g_{13}\theta^{3}u^{3}+g_{14}\theta^{4}u^{4})(v^{p})^{2}  \\
&= (g_{0}+g_{1}\theta u+g_{2}\theta^{2} u^{2}+g_{3}\theta^{3}u^{3}+g_{4}\theta^{4}u^{4}) \\
&\hskip .75cm +(g_{5}\beta+g_{6}\theta\beta u+g_{7}\theta^{2}\beta u^{2}+g_{8}\theta^{3}\beta u^{3}+g_{9}\theta^{4}\beta u^{4})v \\
& \hskip .75cm +(g_{10}\beta^{2}+g_{11}\theta\beta^{2} u+g_{12}\theta^{2}\beta^{2} u^{2}+g_{13}\theta^{3}\beta^{2}u^{3}+g_{14}\theta^{4}\beta^{2}u^{4})v^{2}.
\end{align*}
We precomputed these following values:
$c_{0} = \theta^{2}$
$c_{1} = \theta^{3}$,
$c_{2}=  \theta^{4}$,
$c_{3}= \beta^{2}$,
$c_{4}= \theta \beta$,
$c_{5}= c_{0}\beta$,
$c_{6}= c_{1} \beta$,
$c_{7}= c_{2} \beta $,
$ c_{8}= \theta c_{3}$,
$c_{9}=c_{0} c_{3} $,
$c_{10}=c_{1} c_{3}$,
$c_{11}=c_{2} c_{3}$.
So
\begin{align*}
\pi(a) &=(g_{0}+g_{1}\theta u+g_{2}c_{0} u^{2}+g_{3}c_{1}u^{3}+g_{4}c_{2}u^{4}) \\
&\hskip .75cm +(g_{5}\beta+g_{6}c_{4} u+g_{7}c_{5} u^{2}+g_{8}c_{6} u^{3}
+g_{9}c_{7} u^{4})v \\
&\hskip .75cm +(g_{10}c_{3}+g_{11}c_{8} u+g_{12}c_{9} u^{2}+g_{13}c_{10}u^{3}+g_{14}c_{11}u^{4})v^{2}.
\end{align*}
%The cost of $p$-Frobenius: $14m_{576}+12a_{576}$.\\
%
%\caption{Table 3: Costs of arithmetic operations in $\mathbb{F}_{p^{5}}$ with Newton interpolation and Itoh-Tsuji method}
%
%Let $a \in \mathbb{F}_{p^{15}}$;  $a= a_{0}+ a_{1}v+ a_{2}v^{2}$ with $ a_{i}\in\mathbb{F}_{p^{5}}.$\\
%
The cost of $p$-Frobenius is $14M_1$. This is the same cost as computing $p^2,p^3,p^4$, $p^6,p^7,p^8,p^9$-Frobenius.
For the $p^5$-Frobenius operator, observe from \ref{ci15} that  $v^{p^{5}}=\omega v$ so that
\begin{align*}
a^{p^{5}} &=  (g_{0}+g_{1}u+g_{2}u^{2}+g_{3}u^{3}+g_{4}u^{4})+ (g_{5}+g_{6}u+g_{7}u^{2}+g_{8}u^{3}+g_{9}u^{4})v^{p^{5}} \\
&\hskip .75cm + (g_{10}+g_{11}u +g_{12}u^{2}+g_{13}u^{3}+g_{14}u^{4})(v^{p^{5}})^{2} \\
&= (g_{0}+g_{1}u+g_{2}u^{2}+g_{3}u^{3}+g_{4}u^{4})+ (g_{5}\omega + g_{6}\omega u+g_{7}\omega u^{2}+g_{8}\omega u^{3}+g_{9}\omega u^{4})v \\
&\hskip .75cm + (g_{10}\omega^{2}+g_{11}\omega^{2}u+g_{12}\omega^{2}u^{2}+g_{13}\omega^{2}u^{3}+g_{14}\omega^{2}u^{4})v^{2}.
\end{align*}
We precomputed $d=\omega^{2}$ and
\begin{align*}
\pi^{5}(a) &= a^{p^{5}} \\
&= (g_{0}+g_{1}u+g_{2}u^{2}+g_{3}u^{3}+g_{4}u^{4})+ (g_{5}\omega + g_{6}\omega u+g_{7}\omega u^{2}+g_{8}\omega u^{3}+g_{9}\omega u^{4})v \\
&\hskip .75cm +(g_{10}d+g_{11}du+g_{12}du^{2}+g_{13}du^{3}+g_{14}du^{4})v^{2}.
\end{align*}
The cost of $p^{5}$-Frobenius is $10M_{1}$. This is the same as the cost of $p^{10}$-Frobenius.

%%%%%%%%%%%%%%%%%%%%%%
%%%%%%%%%%%%%%%%%%%%%%
%%%%%%%%%%%%%%%%%%%%%%
%%%%%%%%%%%%%%%%%%%%%%
\end{document}